%
%
%

\input amstex
\documentstyle{amsppt}
\PSAMSFonts

\pageheight{18.5cm}
\magnification=\magstep1
\frenchspacing


\loadbold

\font\smallbf=cmb10 at 8pt
\font\small=cmr10 at 8pt

\def\mapdown#1{\downarrow\rlap{$\vcenter{\hbox{$\scriptstyle#1$}}$}}

\def\phi{\varphi}

\def\nin{\newline\indent}
\def\name#1{{\smc #1\/}}
\def\CC{{\Bbb C}}
\def\FF{{\Bbb F}}

\def\PP{{\Bbb P}}

\def\ZZ{{\Bbb Z}}

\def\EEE{{\Cal E}}

\def\RRR{{\Cal R}}
\def\OOO{{\Cal O}}
\def\QQQ{{\Cal Q}}
\def\VVV{{\Cal V}}

\def\SSS{{\Cal S}}

\def\TTT{{\Cal T}}

\def\KKK{{\Cal K}}

\def\qqed{{\hfill\hfill\qed}}

\def\se#1,#2,#3;{(ab)^{#1}\, a_x^{\,#2}b_x^{\,#3}}

\def\see#1,#2,#3;{(ab)^{#3} (bc)^{#1} (ca)^{#2}\,
\ell_a^{p-#2-#3}\ell_b^{q-#1-#3}\ell_c^{r-#1-#2}}

\def\seep#1,#2,#3;{(ab)^{#3} (bc)^{#1} (ca)^{#2}\,
\ell_a^{p'-#2-#3}\ell_b^{q'-#1-#3}\ell_c^{r'-#1-#2}}

\def\seed#1,#2,#3;{(ab)^{#3} (bc)^{#1} (ca)^{#2}\,
a_x^{p_a-#2-#3}b_x^{p_b-#1-#3}c_x^{p_c-#1-#2}}

\def\SSS{{\Cal S}}
\def\nse#1,#2,#3,#4,#5,#6;{(#1 #2)(#3 #4){#5}_x{#6}_x}
\def\nsee#1,#2,#3,#4,#5,#6;{(#1 #2)(#3 #4)#5_x#6_x}
\def\nsse#1,#2,#3,#4,#5,#6,#7,#8;{(#1 #2)(#3 #4)(#5 #6)(#7 #8)}

\def\tv#1,#2,#3;{[#1,#2]_{#3}}

\TagsOnRight

\input xy
\xyoption{all}

\topmatter
\title Geometry of orbit closures for the representations associated to gradings of Lie algebras of types $E_6$, $F_4$ and $G_2$.
\endtitle
\date December 04, 2016 \enddate
\author Witold Kra\'skiewicz and Jerzy Weyman \endauthor
\address Nicholas Copernicus University
\nin Toru\'n , Poland
\endaddress
\email wkras\@mat.umk.pl \endemail
\address Department of Mathematics, Northeastern University
\nin 360 Huntington Avenue,  BOSTON,  MA 02115, USA \endaddress
\email j.weyman\@neu.edu \endemail

\thanks The study was partially supported by research fellowship within the project “Enhancing Educational Potential of Nicolaus Copernicus University in the Disciplines of Mathematical and Natural Sciences”   (Project № POKL.04.01.01-00-081/10).
 The  second author was partially supported by NSF grant DMS-0600229
\endthanks
\rightheadtext{Geometry of orbit closures for $E_6$, $F_4$, $G_2$}
\abstract In this paper we  investigate the orbit closures for the class of representations of simple algebraic groups
associated to various gradings on a simple Lie algebras of type $E_6$, $F_4$ and $G_2$. The methods for classifying the
orbits for these actions were developed by Vinberg \cite{V75}, \cite{V87}.
We give the orbit descriptions, the degeneration partial orders, and decide the normality of the orbit closures.
We also investigate the rational singularities, Cohen-Macaulay and Gorenstein properties for the orbit closures.
We give the generators of the defining ideals of orbit closures.

\endabstract
\endtopmatter

\document
\def\iz#1:#2.{\noindent\S\ignorespaces#1\dotfill#2\par}

\head  Introduction \endhead

The irreducible representations of reductive groups with finitely many orbits were classified by Kac in \cite{K82} (with some corrections in \cite{DK85}).
They correspond (with very few exceptions) to certain gradings on the root systems, and to the corresponding $\theta$ groups.
In Kac's paper the list of these representations appears in Table II (there are the other tables III, IV classifying so-called visible representations).
We refer to these representations as representations of type I.

 The representations of type I are parametrized by the pairs $(X_n, x)$ where $X_n$ is a  Dynkin diagram with a distinguished node
 $x\in X_n$. This data defines a grading
 $${\goth g}=\oplus_{i=-s}^s {\goth g}_i$$
 of a simple algebra $\goth g$ of type $X_n$ such that the Cartan subalgebra $\goth h$ is contained in ${\goth g}_0$ and the root space ${\goth g}_\beta$ is contained in ${\goth g}_i$ where $i$ is the coefficient of the simple root $\alpha$ corresponding to the node $x$  in the expression for $\beta$ as a linear combination of simple roots. The representation corresponding to $(X_n ,x)$ is the ${\goth g}_1$ with the action of the group $G_0\times \CC^*$ where $G_0$ is the adjoint group corresponding to ${\goth g}_0$ and $\CC^*$ is the copy of $\CC^*$ that occurs in maximal torus of $G$ (the adjoint group corresponding to $\goth g$)  but not in maximal torus of $G_0$.

The orbit closures for the representations of type I were described in two ways by Vinberg in \cite{V75}, \cite{V87}.
The first description states that the orbits are the irreducible components of the intersections of the nilpotent orbits in $\goth g$
with the graded piece ${\goth g}_1$.
In the second paper Vinberg gave a more precise description in terms of the support subalgebras which are graded Lie subalgebras
of the graded Lie algebra $\goth g =\oplus {\goth g}_i$

In this paper we  concentrate on the cases when $X_n$ is equal to $E_6$, $F_4$ and $G_2$.

Besides listing the orbits we give  the degeneration partial order on the orbit closures and investigate
the properties of singularities of the orbit closures.

It turns out that every orbit closure in our representations has a desingularization which is a total space of a homogeneous vector bundle
on some homogeneous space $G_0/P$. This means we can use the geometric method  of \cite{W03}.

The main result of the paper is the calculation of equivariant Euler characteristics of the exterior powers of the bundle
and the corresponding Koszul complex. This gives the formulas for the degeneracy loci of the corresponding orbit closures
in equivariant $K$-theory and gives an indication of what should be the terms of finite free resolution of the coordinate ring
of an orbit closure treated as a module over the coordinate ring of our representation.

By forgetting the equivariant structure we get the calculation of the Hilbert polynomials of the normalizations of the coordinate rings
of the orbit closures.
We are also able to decide the normality, Cohen-Macaulay and Gorenstein
properties of the orbit closures as well as rational singularities property.

Using the results of this paper Federico Galetto \cite{G11} was able to describe explicitly the free resolutions of the coordinate rings of orbit closures as modules over the
coordinate ring of the representation itself.  We list the terms of these resolutions in the corresponding sections. For nonnormal orbit closures $\overline{\OOO_j}$, denoting $N(\overline{\OOO_j})$ the normalization of the orbit closure $\OOO_j$,  we also give the
descriptions of Betti tables of the Cokernel modules $C_j:=\CC[N(\overline{\OOO_j})]/\CC[\overline{\OOO_j}]$. These results were also checked by Macaulay 2 (see \cite{G11}).
See the final remark of section 3 for the outline of the method used to calculate the resolutions. This method gives the generators of the defining ideals in all cases.

The paper is organized as follows. In section 1 we introduce the necessary notation. In section 2 we expose the Vinberg method from \cite {V87}.
Next we give a section surveying the geometric method for calculating syzygies.
In section 4 we deal with the case of degenerate orbits showing how for a degenerate orbit closure the information about syzygies can be deduced from
the information for the smaller orbit closure.
In the following sections we give the explicit descriptions type by type.

 The data are organized as follows. For each representation we start with several  tables.
 First there is a general table with the number of the orbit, the type of the support algebra $\goth s$ and the dimension of the orbit.
 This is followed by  tables with the geometric description of the orbits.
 The numerical data table includes the degree and the numerator of the Hilbert series of the coordinate ring of the normalization of the orbit closure. The denominator is always $(1-t)^{codim}$.
The final table indicates the singularities data i.e. the information on whether the orbit closure and its normalization is spherical, normal, Cohen-Macaulay, has rational singularities and is Gorenstein.

The orbits were calculated first by hand but then the calculations were checked using the program \cite{dG11} kindly provided by Willem de Graaf.
The dimensions of the orbits were calculated using computer routines written by Jason Ribeiro. Jason also wrote a very useful
 python package \cite{R10} for calculating the Euler characteristics of bundles involved.

 The bulk of the calculations was done using  more complicated roots and weight programs (written by the first author) which searched through all parabolic subgroup submodules
 for the representations in question. Then the program calculated needed Euler characteristics and Hilbert polynomials.

In the final section we list some general conclusions about the orbit closures in the representations we deal with.
The corresponding results for the gradings of Lie algebras of types $E_7$ and $E_8$ will appear in subsequent papers.

\bigskip\noindent
{\smallbf Acknowledgment:} {\small
Both authors thank Federico Galetto, Willem de Graaf and Jason Ribeiro
for very useful computer programs.
The second author would like to thank Joe Landsberg and Steven V. Sam
 for interesting conversations and pointing out various errors in earlier versions of this paper.}

\bigskip\bigskip

\head  \S 1. The representations of type I and $\theta$ groups \endhead

Let $X_n$ be a Dynkin diagram and let $\goth g$ be the corresponding
 simple Lie algebra. Let us distinguish a node $x\in X_n$. Let $\alpha_k$ be
a corresponding simple root in the root system $\Phi$ corresponding
to $X_n$. The choice of $\alpha_k$ determines a $\ZZ$- grading on $\Phi$
by letting the degree of a root $\beta$ be equal to the coefficient of
$\alpha_k$ when we write $\beta$ as a linear combination of simple roots. On the level
of Lie algebras this corresponds to a $\ZZ$-grading
$${\goth g}=\oplus_{i\in \ZZ}\ {\goth g}_i .$$

We define the group $G_0 := (G ,G )\times {\CC}^*$ where $(G,G)$ is a
connected semisimple group with the Dynkin diagram $X_n\setminus x$.
A representation of type I is the representation of $G_0$ on ${\goth g}_1$.

We will denote the representation ${\goth g}_1$ by $(X_n,\alpha_k)$.

Denoting by ${\goth l}$ the Levi factor ${\goth g}_0$ we have
$${\goth l}={\goth l}'\oplus {\goth z}(\goth l )$$
where ${\goth l}'$ denotes the Lie algebra associated to $X_n$ with the omitted node $x$, and ${\goth z}({\goth l})$ is a one dimensional center of $\goth l$.

 \name{Vinberg} in \cite{V75} gave a description of the $G_0$-orbits in the representations of type $I$  in terms of
conjugacy classes of nilpotent elements in $\goth g$. Let $e\in {\goth g}_1$
 be a nilpotent element in ${\goth g}$. Consider the irreducible components
of the intersection of the conjugacy class of $e$ in $\goth g$
$$C(e)\cap {\goth g}_1 = C_1 (e)\cup\ldots\cup C_{n(e)}(e) .$$
The sets $C_i (e)$ are clearly $G_0$-stable.  \name{Vinberg's} result
shows that these are precisely the $G_0$-orbits in ${\goth g}_1$.

\proclaim {Theorem 1.1 } The $G_0$-orbits of the action of $G_0$ on
${\goth g}_1$ are the components $C_i (e)$, for all choices of the conjugacy
classes $C(e)$ and all $i$, $1\le i\le n(e)$.
\endproclaim

Theorem 1.1. makes a connection between the orbits in ${\goth g}_1$ and the nilpotent orbits in $\goth g$.
The classification of nilpotent orbits in simple Lie algebras was obtained by Bala and Carter in the papers \cite{BC76a}, \cite{BC76b}.
A good account of this theory is the book \cite{CM93}. Here we recall that the nilpotent orbit of an element $e$ in a simple Lie algebra $\goth g$ is characterized by the smallest Levi subalgebra $\goth l$ containing $e$. One must be careful because sometimes $\goth l$ is equal to $\goth g$. If the element $e$ is a principal element in $\goth l$, then this orbit is denoted by the Dynkin diagram of $\goth l$ (but there might be different ways in which the root system $R({\goth l })$ sits as a subroot system of $R({\goth g})$).

There are, however the nonprincipal nilpotent orbits that are not contained in a smaller reductive Lie algebra $\goth l$. These are called {\it the distinguished nilpotent orbits} and are described in sections 8.2, 8.3, 8.4 of \cite{CM93}. They are characterized by their associated parabolic subgroups (as their Dynkin characteristics are even, see section 8 in \cite{CM93}).
Let us remark that for Lie algebras of classical types, for type $A_n$ the only distinguished nilpotent orbits are the principal ones, and for types $B_n$, $C_n$, $D_n$ these are orbits corresponding to the partitions with different parts. For exceptional Lie algebras the distinguished orbits can be read off the tables in section 8.4. of \cite{CM93}.

Theorem 1.1 is not easy to use because it is not very explicit. In the next section we describe more precise method from second Vinberg paper \cite{V87}.

\bigskip\bigskip

\head \S 2. The Vinberg method for classifying orbits. \endhead

In this section we describe the second paper of Vinberg \cite{V87} in which he describes orbits of nilpotent elements in ${\goth g}_1$.
Similarly to Bala-Carter classification, the nilpotent elements in ${\goth g}_1$ are described by means of some graded subalgebras of $\goth g$.
We need some preliminary notions.

All Lie algebras $\goth g$ we will consider will be Lie algebras of some algebraic group $G$.

Let $(X_n ,\alpha_k )$ be one of the representations on our list.
It defines the grading
$${\goth g}=\oplus_{i\in{\ZZ}} {\goth g}_i$$
where ${\goth g}_i$ is the span of the roots which, written as a combination of simple roots, have $\alpha_k$ with coefficient $i$.
The component ${\goth g}_0$ contains in addition a Cartan subalgebra. $G_0$ denotes the connected component of the subgroup of $G$ corresponding to ${\goth g}_0$.

In the sequel $Z(x)$ denotes the centalizer of an element $x\in G$. $Z_0(x)$ is $Z(x)\cap G_0$. The gothic letters, ${\goth z}$, ${\goth z}_0$ denote corresponding Lie algebras.
Similarly, $N(x)$ denotes the normalizer of an element $x\in G$. $N_0(x)$ is $N(x)\cap G_0$.

We denote $R({\goth g})$ the set of roots of a reductive Lie algebra $\goth g$, and $\Pi ({\goth g})$ denotes the set of simple roots.

In order to state Vinberg Theorem we need some definitions.

We will be dealing with the graded Lie subalgebras ${\goth s}=\oplus_{i\in\ZZ}{\goth s}_i$.

\proclaim {Definition 2.1} A graded Lie subalgebra $\goth s$ of $\goth g$ is regular if it is normalized by a maximal torus in ${\goth g}_0$.
A reductive graded Lie algebra $\goth s$ of $\goth g$ is complete if it is not  a proper graded Lie subalgebra of any regular reductive $\ZZ$-graded Lie algebra of the same rank.
\endproclaim

\proclaim {Definition 2.2} A $\ZZ$-graded Lie algebra $\goth g$ is locally flat if one of the following equivalent conditions is satisfied, for $e$ a point in general position in ${\goth g}_1$:
\item{(1)} The subgroup $Z_0(e)$ is finite,
\item{(2)} ${\goth z}_0(e)=0$,
\item{(3)} $dim\ {\goth g}_0 = dim\ {\goth g}_1$.
\endproclaim

Let ${\goth g}=\oplus_{i\in\ZZ}{\goth g}_i$ be a graded reductive Lie algebra. Let us fix a nonzero nilpotent element $e\in {\goth g}_a$. Let us choose some maximal torus $H$ in $N_0(e)$.
Its Lie algebra $\goth h$ is {\it the accompanying torus} of the element $e$. We denote by $\phi$ the character of the torus $H$ defined by the condition
$$[u,e]=\phi(u)e$$
for $u\in{\goth h}$. Consider the graded Lie subalgebra ${\goth g}({\goth h}, \phi )$ of $\goth g$ defined as follows
$${\goth g}({\goth h}, \phi ) =\oplus_{i\in\ZZ}\ {\goth g}({\goth h}, \phi )_i$$
where
$${\goth g}({\goth h}, \phi )_i=\lbrace x\in {\goth g}_{ia} \ |\ [u,x]=i\phi(u)x\ \forall u\in H\ \rbrace .$$

\proclaim {Definition 2.3} The support $\goth s$ of the nilpotent element $e\in {\goth g}_a$ is the commutant of ${\goth g}({\goth h}, \phi )$ considered as a $\ZZ$ graded Lie algebra.
\endproclaim

Clearly $e\in {\goth s}_1$.

We are ready to state the main theorem of \cite {V87}.

\proclaim {Theorem 2.4. (Vinberg)} The supports of nilpotent elements of the space ${\goth g}_i$ are exactly the complete regular locally flat semisimple $\ZZ$-graded subalgebras of the algebra $\goth g$.
The nilpotent element $e$ can be recovered from the support subalgebra $\goth s$ as the generic element in ${\goth s}_1$.
\endproclaim

 It follows from the theorem that the nilpotent element $e$ is defined uniquely (up to conjugation by an element of $G_0$) by its support.
 This means it is enough to classify the regular semisimple $\ZZ$-graded subalgebras $\goth s$ of $\goth g$.

 Let us choose the maximal torus $\goth t$ of ${\goth g}_0$. The $\ZZ$-graded subalgebra $\goth s$ is {\it standard } if it is normalized by $\goth t$, i.e. if for all $i\in\ZZ$ we have

 $$[{\goth t}, {\goth s}_i]\subset {\goth s}_i .$$

Vinberg also proves that every $\ZZ$ graded subalgebra $\goth s$ is conjugated to a standard subalgebra by an element of $G_0$. Moreover, he shows that if two standard
$\ZZ$-graded subalgebras are conjugated by an element of $G_0$, then they are conjugated by an element of $N_0 ({\goth t})$.

This allows to give a combinatorial method for classifying regular semisimple $\ZZ$-graded subalgebras of $\goth g$.

Let $\goth s$ be a standard semisimple $\ZZ$-graded subalgebra of $\goth g$. The subalgebra $\goth s$ defines the degree map $deg: R({\goth s)}\rightarrow\ZZ$.
For a standard $\ZZ$-graded subalgebra $\goth s$ we also get the map
$$f: R({\goth s)}\rightarrow R({\goth g}).$$

The map $f$ has to be {\it additive}, i.e. it satisfies

$$f(\alpha +\beta )=f(\alpha )+f(\beta )\ \ \forall\alpha ,\beta\in R({\goth s}),$$
$$f(-\alpha  )=-f(\alpha )\ \ \forall\alpha \in R({\goth s}).$$

Moreover we have

\proclaim {Proposition 2.5} The map $f$ satisfies the following properties:
\item{a)}
$$ {{(f(\alpha ),f(\beta ))}\over{(f(\alpha ),f(\alpha ))}}= {{(\alpha ,\beta )}\over{(\alpha ,\alpha )}}\ \ \forall\alpha, \beta\in R({\goth s}).$$
\item{b)} $f(\alpha )-f(\beta )\notin R({\goth g})\ \ \forall\alpha ,\beta\in \Pi({\goth s})$,
\item{c)} $deg\ f(\alpha )=deg\ \alpha,\ \ \forall\alpha\in \Pi({\goth s})$.

Conversely, every map satisfying $a)$, $b)$, $c)$ defines a standard regular $\ZZ$-graded subalgebra ${\goth s}$ of ${\goth g}$.
\endproclaim

\proclaim {Remark 2.6} The subalgebra $\goth s$ corresponding to the map $f$ is complete if and only if there exists an element $w$ in the Weyl group $W$  of $\goth g$
such that $wf(\Pi({\goth s}))\subset\Pi({\goth g})$ (see \cite{V87}, p.25) .
\endproclaim

Proposition 2.5 means that in order to classify the nilpotent elements $e\in {\goth g}_1$ we need to classify the possible maps $f$
corresponding to its support, i.e. the corresponding complete regular $\ZZ$-graded subalgebra $\goth s$.
Since we are interested in the nilpotents $e\in {\goth g}_1$, we need to classify the maps $f$ for which $deg( f(\alpha ))\in\lbrace 0,1\rbrace$ for every $\alpha\in \Pi({\goth s})$.

At that point we can make the connection with the Bala-Carter classification more explicit. The Bala-Carter characteristic of a nilpotent element is given by the type of the
algebra $\goth s$. So in the case of representations of type I it follows that in order to find the components of a given nilpotent orbit with Bala-Carter characteristic $\goth s$ it is enough to see in how many ways
the Lie algebra $\goth s$ can be embedded as a $\ZZ$-graded Lie algebra into $\goth g$.

\bigskip\bigskip

In each case below we follow the following strategy.

\proclaim {Proposition 2.7}
Let $\goth g$ be a simple Lie algebra. Let ${\goth g}=\oplus_{i=-s}^s {\goth g}_i$ be the grading associated to a given simple root.
The following procedure allows to find the components of the intersction of a nilpotent orbit in $\goth g$ with ${\goth g}_1$.
\item{1.} Establish the restriction of the invariant scalar product $( ,)$ on $\goth g$ to ${\goth g}_1$.
\item{2.} Fix a nilpotent element $e\in{\goth g}$ with Bala-Carter characteristic $\goth s$.
\item{3.} Find  in how many ways  $\goth s$ embeds into $\goth g$ as a standard $\ZZ$-graded Lie subalgebra by exhibiting corresponding map $f$
as in Proposition 2.5.
\item{4.} For a characteristic $\goth s$ such that $e$ is a principal nilpotent in $\goth s$  the map $f$ sends all simple roots of $\goth s$ to weight vectors in ${\goth g}_1$.
Thus it is enough to classify the subsets of weights vectors in ${\goth g}_1$ for which the pattern of scalar products is the same as the one for simple roots of $\goth s$.
\item{5.} For non-principal orbits one needs to make a more detailed analysis, but still one reduces to finding sets of weight vectors in ${\goth g}_1$ with certain patterns of scalar products.
\endproclaim

\bigskip\bigskip

\head \S 3. The geometric technique. \endhead

In this section we provide a quick description of main facts related to geometric technique of calculating syzygies.
We work over an algebraically closed field $K$.

Let us consider the projective variety $V$ of dimension $m$.
Let $X=A^N_\CC$ be the affine space. The space $X\times V$ can be viewed
as a total space of trivial vector bundle $\EEE$ of rank $N$ over $V$.
Let us consider the subvariety $Z$ in $X\times V$ which is the total space of a subbundle
$\SSS$ in $\EEE$.
We denote by $q$ the projection $q: X\times V\longrightarrow X$ and by $q^\prime$ the restriction
of $q$ to $Z$. Let $Y=q(Z)$. We get the basic diagram
$$                  \matrix Z&\subset&X\times V\\
                 \mapdown{q^\prime}&&\mapdown{q}\\
                 Y&\subset&X
                 \endmatrix$$

The projection from $X\times V$ onto $V$ is denoted by $p$ and the quotient bundle ${\EEE}/{\SSS}$
by $\TTT$.
Thus we have the exact sequence of vector bundles on $V$
$$           0\longrightarrow {\SSS}\longrightarrow {\EEE}\longrightarrow {\TTT}\longrightarrow 0.$$
 The ranks of $\SSS$ and $\TTT$ will be denoted by $s$, $t$ respectively.
The coordinate ring of $X$ will be denoted by $A$. It is a polynomial ring in
 $N$ variables over $\CC$.
We will identify the sheaves on $X$ with $A$-modules.

 The locally free resolution of the sheaf ${\OOO}_Z$ as an ${\OOO}_{X\times V}$-module is given by the Koszul complex
$$ {\KKK}_{\bullet} (\xi ): 0\rightarrow \bigwedge^t (p^* \xi )\rightarrow \ldots
\rightarrow \bigwedge^2 (p^* \xi )\rightarrow p^* (\xi ) \rightarrow {\OOO}_{X\times V}$$
where $\xi = {\TTT}^*$. The differentials in this complex are homogeneous of degree $1$ in
the coordinate functions on $X$.
  The direct image $p_* ({\OOO}_Z )$ can be identified with the
the sheaf of algebras $Sym(\eta )$ where $\eta = {\SSS}^*$.

The idea of the geometric technique is to use the Koszul complex ${\KKK}(\xi )_{\bullet} $ to
construct for each vector bundle $\VVV$ on $V$ the free complex $\FF_{\bullet} ({\VVV})$ of
$A$-modules with the homology supported in $Y$. In many cases  the complex ${\FF}({\OOO}_V )_{\bullet} $
gives the free resolution of the defining ideal of $Y$.

For every vector bundle $\VVV$ on $V$ we introduce the complex
$${\KKK}(\xi ,{\VVV})_{\bullet}  := {\KKK}(\xi )_{\bullet} \otimes_{{\OOO}_{X\times V}} p^* {\VVV}$$ This complex is a locally free resolution of the ${\OOO}_{X\times V}$-module
$ M({\VVV}) :={\OOO}_Z \otimes p^* {\VVV}$.

Now we are ready to state the basic theorem (Theorem (5.1.2) in [W]).

\proclaim  {Theorem 3.1}  For a vector bundle $\VVV$ on $V$
we define a free graded $A$-modules
$$ \ {\FF}({\VVV})_i  = \bigoplus_{j\ge 0} H^j (V,\bigwedge^{i+j}\xi\otimes{\VVV}
)\otimes_{k} A(-i-j)$$
 \item {a)} There exist  minimal differentials
$$ d_i ({\VVV}): {\FF}({\VVV})_i \rightarrow {\FF}({\VVV})_{i-1} $$
of degree $0$ such that ${\FF}({\VVV})_{\bullet} $ is a complex of graded free
$A$-modules with
$$ H_{-i} ({\FF}({\VVV})_{\bullet} ) = {\RRR}^i q_* M({\VVV})$$
In particular the complex ${\FF}({\VVV})_{\bullet} $ is exact in positive degrees.
\item {b)} The sheaf ${\RRR}^i q_* M({\VVV})$ is equal to $H^i (Z, M({\VVV}))$ and it
can be also identified with the graded $A$-module
$H^i (V, Sym (\eta )\otimes {\VVV})$.
\item {c)} If $\phi : M({\VVV})\rightarrow M({\VVV}^\prime )(n)$ is a morphism of
graded sheaves then there exists a morphism of complexes
$$f_{\bullet} (\phi ): {\FF}({\VVV})_{\bullet} \rightarrow {\FF}({\VVV}^\prime )_{\bullet} (n)$$
Its induced map $H_{-i} (f_{\bullet} (\phi ))$ can be identified with the induced map
$$H^i (Z, M({\VVV}) )\rightarrow H^i (Z, M({\VVV}^\prime ) )(n).$$
\endproclaim

\vskip .2cm

If $\VVV$ is a one dimensional trivial bundle on $V$ then the complex ${\FF}({\VVV})_{\bullet} $
is denoted simply by ${\FF}_{\bullet}$.

  The next theorem gives the criterion for the complex ${\FF}_{\bullet}$ to be the free resolution of the
coordinate ring of $Y$.

\proclaim  {Theorem 3.2} Let us assume that the map $q^\prime : Z\longrightarrow Y$
is a birational isomorphism. Then the following properties hold.
\item {a)} The module $q^\prime_* {\OOO}_Z$ is the normalization of ${\CC}[Y]$.
\item {b)} If ${\RRR}^i q^\prime_* {\OOO}_Z = 0$ for $i>0$, then ${\FF}_{\bullet}$ is a finite
free resolution of the normalization of ${\CC}[Y]$ treated as an $A$-module.
\item {c)} If ${\RRR}^i q^\prime_* {\OOO}_Z = 0$ for $i>0$ and ${\FF}_0 = H^0 (V, \bigwedge^0 \xi
)\otimes A = A$ then $Y$ is normal and it has rational singularities.
\endproclaim

This is Theorem (5.1.3) in [W].

The complexes $\FF({\VVV})_{\bullet} $ satisfy the Grothendieck type duality. Let $\omega_V$ denote the
canonical divisor on $V$.

\proclaim {Theorem 3.3} Let ${\VVV}$ be a vector bundle on $V$. Let us introduce the dual bundle
$${\VVV}^\vee = \omega_V \otimes \bigwedge^t \xi^* \otimes {\VVV}^* .$$
Then
$$ \FF({\VVV}^\vee )_{\bullet}  = \FF({\VVV})_{\bullet}^* [m-t]$$
\endproclaim

This is Theorem (5.1.4) in [W].

In some of our applications the projective variety $V$ will be a Grassmannian. To fix the notation, let us work with the
Grassmannian $Grass(r, E)$ of subspaces of dimension $r$ in a vector space $F$ of dimension $n$.
Let $$0\rightarrow {\RRR}\rightarrow E\times Grass(r, E)\rightarrow {\QQQ}\rightarrow 0$$
be a tautological sequence of the vector bundles on $Grass(r, E)$.
The vector bundle
$\xi$ will be a direct sum of the bundles of the form $S_{\lambda_1 ,\ldots ,\lambda_{n-r}}{\QQQ}\otimes S_{\mu_1
,\ldots ,\mu_r}{\RRR}$.  Thus all the exterior powers of $\xi$ will also be the direct sums of such bundles. We will
apply repreatedly the following result to calculate cohomology of vector bundles $S_{\lambda_1 ,\ldots
,\lambda_{n-r}}{\QQQ}\otimes S_{\mu_1 ,\ldots ,\mu_r}{\RRR}$.

\proclaim {Proposition 3.4 (Bott's algorithm)}. The cohomology of the vector bundle $S_{\lambda_1 ,\ldots ,\lambda_{n-r}}{\QQQ}\otimes
S_{\mu_1 ,\ldots ,\mu_r}{\RRR}$
 on $Grass(r, E)$ is calculated
as follows. We look at the weight
$$(\lambda ,\mu ) = (\lambda_1 ,\ldots ,\lambda_{n-r},\mu_1 ,\ldots ,\mu_r )$$
 and add to it $\rho = (n, n-1,\ldots
,1)$.  Then one of two mutually exclusive cases occurs.
\item{a)} The resulting sequence
$$(\lambda ,\mu ) +\rho =(\lambda_1 +n,\ldots \lambda_{n-r} +r+1,\mu_1 +r ,\ldots
,\mu_r +1 )$$
 has repetitions. In such case
$$ H^i (Grass(r, E), S_\lambda {\QQQ}\otimes S_{\mu} {\RRR} )=0$$
for all $i\ge 0$.
\item{b)}  The  sequence $(\lambda ,\mu )+\rho$ has no repetitions. Then there is a unique permutation $w\in \Sigma_n$
that makes this sequence decreasing.  The sequence $\nu =w((\lambda ,\mu ) +\rho )-\rho$ is a non-increasing
sequence. Then the only non-zero cohomology group of the sheaf $S_{\lambda}{\QQQ}\otimes S_{\mu} {\RRR}$ is
 the group $H^l$ where $l=l(w)$ is the length of $w$. We have
 $$H^l (Grass(r, E), S_{\lambda}{\QQQ}\otimes S_{\mu} {\RRR} )=S_\nu E.$$
Here $S_\nu E$ denotes the highest weight representation $S_\nu E$ of $GL(E)$ corresponding
to the highest weight $\nu$ (so-called Schur module).
\endproclaim

\proclaim {Remark 3.5 (desingularizations and confirming the resolution)}
\item{a)} Desingularizations $Z$ of the orbit closures were determined by brute force. We went through the list of all possible $B$-submodules of the representations $V$ and calculated
the equivariant Euler characteristics of the exterior powers of the corresponding bundles $\xi$. This allowed to make a guess as to what the desingularization $Z$ is.
Checking it is not difficult, first one calculates the dimension and then the fiber over the representative of the orbit listed in the table.
We do not give these calculations due to the number of cases.
\item{b)}
In the cases of orbit closures investigated in this paper the vector bundles $\xi$ are not semisimple. This makes the
calculation of cohomology difficult. We calculated the equivariant Euler characteristics of the
exterior powers of the bundles $\xi$. This gives the formula for the equivariant $K$-theory class of the degeneracy locus corresponding to a given orbit closure.
This also gives the representations which have to occur in the resolution.
However one might very well have additional ghost terms which might occur. The ghost terms occur in pairs (in the same homogeneous degree but in neighboring homological degrees)
so they cancel out in the Euler characteristic calculation.
To confirm that such ghost terms do not appear we checked the acyclicity of the complexes in question by different  method.
First we attempted to calculate the resolution by Macaulay 2. Usually the calculation was too big to finish. However by knowing the degree
of the terms expected from the Euler characteristic calculation we could produce the individual matrices in the expected resolutions and prove that they form a complex. To prove the acyclicity of the complexes obtained in this way  we used the Buchsbaum-Eisenbud exactness criterion \cite{BE73}.This criterion states that in order to prove acyclicity one needs to check the ranks of the maps in the complex and the depths of the ideals of the maximal non-vanishing minors of the matrices in our complex.  However we know by construction that our complexes are equivariant. This means that
in order to calculate ranks and depths of the ideals of minors one just needs to calculate the ranks
of the matrices in the resolution specialized to representatives of the orbits. This was done by computer.
Also claims involving the cokernel modules for the nonnormal orbit closures were checked by Macaulay 2 in similar way.
These calculations were performed by Federico Galetto and the description of the results will appear elsewhere.

\endproclaim

\proclaim{Remark 3.6. ( proving normality and rational singularities)}
In the cases we consider below we claim the general fact that the normalizations of orbit closures have rational singularities and we list the normal orbit closures.
These results are established in the following way. For each orbit closure $\overline{\OOO}$ we exhibit its desingularization $Z$ which is the total space of the homogeneous vector bundle
on some homogeneous space $G/P$. This gives us the bundles $\eta$ and $\xi$ in each case.
We have the following vanishing results which allow to claim the normality and rational singularities results. The non-normality of the remaining orbit closures is easy to establish
by finding that the term $\FF_0$ has non-trivial summands in positive degree.
\endproclaim
\proclaim{Proposition 3.7} For each orbit closure $\overline{\OOO}$ and its desingularization $Z$ we have
\item {a)} $H^i(G/P, Sym(\eta))=0$ for $i>0$,
\item {b)} For the orbits that are marked as normal in the tables we have $H^i (G/P, \bigwedge^j\xi)=0$ for all $i\ge j$.
\endproclaim
\demo{Proof} The idea of the proof is as follows. Consider the pull-backs of the bundles $\eta$ and $\xi$ to $G/B$.
We can perform the cohomology calculations there. The bundles $\eta$ and $\xi$ can be filtered by line bundles on $G/B$.
Denote by $\eta'$ and $\xi'$ the associated graded bundles which are the sums of line bundles on $G/B$.
One has the following curious combinatorial fact, checked case by case by computer.
\proclaim{Proposition 3.8}
Let $j>0$.
There is no irreducible representation $V_\lambda$ of $G$ occurring in two  $G$-modules:  $H^i(G/B, \bigwedge^j\xi')$ for $i>j$,  and in $H^* (G/B, S_j\eta')$. If the orbit closure $\overline{\OOO}$ is marked as normal in the tables, there is no irreducible representation $V_\lambda$ of $G$ occurring in the $G$-module $H^i(G/B, \bigwedge^j\xi')$ for $i\ge j$ and also occurring
in the $G$-module $H^* (G/B, S_j\eta')$.
\endproclaim
The Proposition 3.8  implies our results. Indeed, assume that the normalization of $\overline{\OOO}$ has no rational singularities.
Then the complex $\FF_\bullet$ has to have terms in negative degrees. Let $\FF_{-s}$ be the smallest non-zero term of $\FF_\bullet$.
Let $V_\lambda$ be the irreducible representation occurring in $\FF_{-s}$. Then this representation has to occur (for the appropriate $j$) both in
$H^s(G/B, S_j\eta')$ and in $H^{j+s}(G/B, \bigwedge^j\xi')$. This gives the contradiction with the first statement in the Proposition 3.8.
Assume now that the orbit closure marked as normal in the tables is not in fact normal. Then the complex $\FF_0$ has a non-trivial term in positive degree.
Let $V_\lambda$ be the irreducible representation occurring in $\FF_{0}$ in homogeneous degree $j$. Then this representation has to occur both in
$H^0(G/B, S_j\eta')$ and in $H^{j}(G/B, \bigwedge^j\xi')$. This gives the contradiction with the second statement in the Proposition 3.8.
\qqed
\enddemo

\bigskip\bigskip

\head \S 4. The degenerate orbits. The highest weight orbits and hyperdiscriminants. \endhead

\proclaim {Definition 4.1}
Let $(X_n, \alpha_k)$ be a Dynkin diagram with a distinguished node. Let
$${\goth g}=\oplus_{-s\le i\le s} {\goth g}_i$$
be the grading on the Lie algebra $\goth g$ of type $X_n$ corresponding to the simple root $\alpha_k$.
An orbit $\OOO_v$ ($v\in{\goth g}_1$) is degenerate, if there exists a representative $w\in \OOO_v$ and a node $\alpha_l$ ($l\ne k$ ) in $X_n$ such that
$w\in {\goth g}'_1$ where ${\goth g}'$ is a Lie subalgebra of $\goth g$ corresponding to the subdiagram  $X_n\setminus\alpha_l$, i.e. if the corresponding orbit occurs already in the smaller case $(X_n\setminus\alpha_l ,\alpha_k)$. We can always assume that $\alpha_l$ is an extremal node of $X_l$, i.e. it is connected to just one other node in $X_n$.
\endproclaim

Let us assume that the orbit closure for $(X_n ,\alpha_k)$ is degenerate, and the orbit comes from the smaller orbit
$\OOO_v$ from $(X_n\setminus\alpha_l ,\alpha_k)$.
 In our situation, ${\goth g}'_1$ is a $P_l$-submodule of ${\goth g}_1$.
We consider the incidence variety

$$Z(X_n,k,l):={\goth g}'_1\times_{P_l} G_0.$$
where the subscript $_{P_l}$ means we divide by the equivalence relation
$$(v,g)\cong (pv, gp^{-1}).$$
The variety $Z(X_n,k,l)$ has two projections
$$p:Z(X_n,k,l)\rightarrow G_0/P_l ,\ \ \ q: Z(X_n,k,l)\rightarrow {\goth g}_1.$$
The first projection makes $Z(X_n,k,l)$ a total space of the homogeneous vector bundle over $G_0/P_l$, with the fibre ${\goth g}'_1$.
The second projection is given by
$(v,g)\mapsto g^.v$. It projects $Z(X_n,k,l)$ onto the subvariety
$$Y(X_n,k,l)=\lbrace v\in {\goth g}_1\ \ |\ \exists g\in G_0\ g^. v\in {\goth g}'_1\rbrace.$$

Now it is clear we are in the setup of the section 3. Thus we have the bundles $\eta (X_n,k,l)$, $\xi (X_n,k,l)$
the Koszul complex and twisted complexes $\FF (X_n,k,l; V)_\bullet$ for any $L$-module $V$. If $V$ is a trivial $L$-module we talk about the complex
$\FF (X_n,k,l)_\bullet$.

\proclaim {Proposition 4.2} Let $\alpha_l$ be an extremal node in the Dynkin diagram $X_n$.
Let $Y_m$ be the connected component of $X_n\setminus\alpha_k$.
We have the following possibilities
\item{a)} $(Y_m,\alpha_l )= (A_m,\alpha_m)$, the bundles $\eta (X_n,k,l)$, $\xi (X_n,k,l)$ are semi-simple,
\item{b)} $(Y_m,\alpha_l)= (D_m,\alpha_{m-1})$, $(Y_m,\alpha_l )=(D_m ,\alpha_m)$, the bundles $\eta (X_n,k,l)$, $\xi (X_n,k,l)$ are semi-simple,
\item{c)} $(Y_m,\alpha_l)= (D_m,\alpha_{1})$,  the bundles $\eta (X_n,k,l)$, $\xi (X_n,k,l)$ are semi-simple,
\item{d)} $(Y_m,\alpha_l )=(E_6,\alpha_1)$,  the bundle $\eta (X_n,k,l)$ is semi-simple, $\xi (X_n,k,l)$ is the extension of a semi-simple bundle by a one dimensional trivial bundle,
\item{e)} $(Y_m,\alpha_l )=(E_6,\alpha_2)$, the bundle $\eta (X_n,k,l)$ is semi-simple, $\xi (X_n,k,l)$ is the extension of a semi-simple bundle by a one dimensional trivial bundle,
\item{f)} $(Y_m,\alpha_l )=(E_7,\alpha_1)$,
\item{g)} $(Y_m,\alpha_l )=(E_7,\alpha_2)$, the bundle $\xi (X_n,k,l)$ is semi-simple.

The cases f) and g) could appear only in the case when ${\goth g}_1$ is $V(\omega_7, E_7)$ but this case is subregular.
So we conclude that the degenerate orbit closures are manageable.
\endproclaim

\demo{Proof}
Let us consider these cases one by one. Case a) is the usual case of incraesing the dimension of the vector sapace.
We have $\eta=\QQQ$, $\xi =\RRR$ (tensored with the representatons corresponding to the remaining connected components of $X_n\setminus\alpha_k$).

Case b) is similar. We deal with an orthogonal vector space $F$ of even dimension. The homogeneous space $G_0/P_l$ the isotropic Grassmannian (or rather its connected component with tautological sequence
$$0\rightarrow \RRR\rightarrow F\times G_0/P_l\rightarrow \QQQ\rightarrow 0$$
where $\QQQ=\RRR^\vee$. In this case $\xi$ is semi-simple so one has an algorithm to calculate the complexes $\FF(D_m, k,l)_\bullet$ and
$\FF(D_m, k,l;V)_\bullet$.

In the case c) the homogeneous space is the quadratic with the tautological filtration
$$0\subset\RRR\subset\RRR^\vee\subset F\times G_0/P_l.$$
Our representation ${\goth g}_1$  is the spinor representation. Without loss of generality we can assume it is the one corresponding
to the fundamental weight $\omega_m$.  The $P_l$-module corresponding to the bundle $\eta (X_n,k,l)$ is the span of all weights with $1\over 2$ on first coordinate, the $P_l$-module corresponding to the bundle $\xi (X_n,k,l)$ is the span of all weights with ${-1}\over 2$ on the first coordinate. Both bundles are semi-simple.

In the case d) the representation ${\goth g}_1$ is the $27$-dimensional representation of ${\goth g}(E_6)$.
After restricting to $L= Spin(10)\times\CC^*$ the representation $V(\omega_1, E_6)$ decomposes as $V(\omega_4, D_5)\oplus V(\omega_1, D_5)\oplus \CC$. The $P_l$-module corresponding to the bundle $\eta (X_n,k,l)$ is the $L$-module $V(\omega_5, D_6)$. The $P_l$-module corresponding to the bundle  $\xi (X_n,k,l)$ is the extension of $V(\omega_1 ,D_6)$ by a one dimensional  trivial module. So $\xi (X_n,k,l)$ is not semi-simple, but still manageable.

In the case e) the representation ${\goth g}_1$ is the $27$-dimensional representation of ${\goth g}(E_6)$.
After restricting to $L= GL(6)$ the representation $V(\omega_1, E_6)$ decomposes as $\bigwedge^3\CC^6\oplus \CC^6\oplus\CC$. The $P_l$-module corresponding to the bundle $\eta (X_n,k,l)$ is the $L$-module $\bigwedge^3\CC^6$. The $P_l$-module corresponding to the bundle bundle $\xi (X_n,k,l)$ is the extension of $\CC^6$ by a one dimensional  trivial module. So $\xi (X_n,k,l)$ is not semi-simple, but still manageable.

The case f) is the most complicated. The representation ${\goth g}_1$ is the $56$-dimensional representation $V(\omega_7, E_7)$ of ${\goth g}(E_7)$.
After restricting to $L= Spin(12)\times\CC^*$ the representation $V(\omega_1, E_6)$ decomposes as $V(\omega_5, D_6)\oplus V(\omega_1, D_6)\oplus V(\omega_1, D_6)$. The $P_l$-module structure is such that $V(\omega_5, D_6)$ is the middle factor in the filtration.

The case g) is a little easier. The representation ${\goth g}_1$ is the $56$-dimensional representation $V(\omega_7, E_7)$ of ${\goth g}(E_7)$.
After restricting to $L= GL(7)$ the representation $V(\omega_1, E_6)$ decomposes as $\bigwedge^3\CC^7\oplus\bigwedge^5\CC^7$. The $P_l$-module corresponding to the bundle $\eta (X_n,k,l)$ is the $L$-module $\bigwedge^3\CC^7$. The $P_l$-module corresponding to the bundle bundle $\xi (X_n,k,l)$ is the $\bigwedge^5\CC^7$. So $\xi (X_n,k,l)$ is semi-simple.

\qqed
\enddemo

Let us assume that the orbit closure for $(X_n ,\alpha_k)$ is degenerate, of type a)-e), and the orbit comes from the smaller orbit
$\OOO_v$ from $(X_n\setminus\alpha_l ,\alpha_k)$. Assume in fact that we have the free resolution of the coordinate ring
$A'/J$ of $\overline{\OOO_v}$ over the coordinate ring $A':= \CC[{\overline{\OOO_v}}]$. Let this resolution have the terms
$$\FF'_\bullet:0\rightarrow \FF'_s\rightarrow\ldots\rightarrow \FF'_1\rightarrow \FF'_0=A'.$$

The terms $\FF'_i$ and the ring $A'$ have the action of the group $L:=L_l$ which is the Levi factor of the parabolic $P_l$ of $G_0$
associated to the root $\alpha_l$. Moreover, in our situation, ${\goth g}'_1$ is a $P_l$-submodule of ${\goth g}_1$.
Extending this action to $P_l$ we see that the terms $p_* \FF'_i$ are of the form $\oplus V(\lambda, L)\otimes Sym(\eta (X_n,k,l))$
for some dominant integral weights of $L$, and for each of these weights $\lambda$  the sheaf on $G_0/P_l$ associated to $V(\lambda ,L)$ has no higher cohomology.

\proclaim {Proposition 4.3}
Let $\OOO_v$ be a degenerate orbit of type a)-e).
\item{a)} The sheaves $\FF'_i$ have no higher cohomology,
\item{b)} The modules $H^0 (G_0/P_l , p_*(\FF'_i))$ are $A$-modules supported
in $Y(X_n,k,l)$. The minimal free resolution over $A$ of each direct summand of the form $V(\lambda, L)\otimes Sym (\eta (X_n,k,l))$ can be calculated using the twisted Koszul complexes
$V(\lambda ,L)\otimes\bigwedge^\bullet (\xi (X_n,k,l))$. By Proposition 4.2 their terms are computable.
\item{c)} The minimal free resolutions of the modules $H^0 (G_0/P_l , p_*(\FF'_i))$ can (by iterated cone construction) be assembled into the terms of a non-minimal free resolution of $\CC[\overline{\OOO_v}]$.
\endproclaim

\demo{Proof}

The terms $p_* \FF'_i$ are of the form $\oplus V(\lambda, L)\otimes Sym(\eta (X_n,k,l))$
for some dominant integral weights of $L$, and for each of these weights $\lambda$  the sheaf on $G_0/P_l$ associated to $V(\lambda ,L)$ has no higher cohomology. Therefore higher cohomology of $\oplus V(\lambda, L)\otimes Sym(\eta (X_n,k,l))$ vanishes.
The rest follows from geometric technique and iterated cone construction.
\qqed
\enddemo

Proposition 4.3 gives a method to find a non-minimal resolution of the coordinate rings of the closures of the degenerate orbit closures. This allows to decide normality, rational singularities, Cohen-Macaulay and Gorenstein properties of degenerate orbit closures.
Thus in dealing with minimal free resolutions of the coordinate rings we will concentrate on non-degenerate orbit closures.

\bigskip

We close out this section with remarks on some typical orbit closures which occur in all cases.
The smallest non-zero orbit is the orbit of the highest weight vector. Geometrically this is a cone over some homogeneous space for the group $G_0$.
The orbits of highest weight vectors were extensively studied. Their main properties are listed in the following proposition.

\proclaim {Proposition 4.4}
Let $V_\lambda$ be an irreducible representation of highest weight $\lambda$ for the reductive group $G$.
Let $X\subset V_\lambda$ be the closure of the orbit $G^.v_\lambda$ where $v_\lambda$ is the highest weight vector in $V_\lambda$.
\item{a)} $K[X]=\oplus_{n\ge 0} V_{n\lambda}^*$,
\item{b)} $X$ is normal, Cohen-Macaulay and has rational singularities,
\item{c)} The defining ideal of $X$ is generated by quadrics.
\endproclaim

\demo{Proof} For a) and b) see \cite{W03}, exercise 8, chapter 5. For the proof of c) (which is known as Kostant Theorem) see the very useful book by Timashev \cite{T06}, section 28.
\qqed
\enddemo

\proclaim {Remark 4.5} In all representations of type I the highest weight orbit  is the orbit $\OOO_1$.
\endproclaim

Another orbit closure which occurs in all cases is {\it  the hyperdiscriminant orbit closure}. This is a subvariety which is projectively dual to $\overline{\OOO_1}$.
It always has a desingularization given by a vector bundle.  For the general properties of such orbits see \cite{W03} chapter 9 and \cite{WZ96}.

Let $G$, $G_0$, ${\goth g}=\oplus_{i=-m}^m {\goth g}_i$ have the same meaning as above.

The bundle $\xi_{disc.}$ is the bundle of $1$-jets. It is spanned (as a $B_0$-module) by the lowest weight vector $v'_\lambda\in{\goth g}^*_1$ and the elements $\alpha v'_\lambda$ of the roots $\alpha$ of the Lie algebra ${\goth g}_0$  on the lowest weight vector $v'_\lambda\in{\goth g}^*_1$.

We can then construct, as in section 3, the vector bundle $Z_{disc.}$ defined by the natural cosection of $\xi_{disc.}$ induced by its embedding into ${\goth g}^*_1$.  We denote $X_{disc.}$ the image $q(Z_{disc.})$.

\bigskip\bigskip

\head \S 5. The type I $E_6$ representations. \endhead

Here we just enumerate all possibilities for possible simple roots $\alpha_k$. We number the nodes of Dynkin diagrams according to the conventions of Bourbaki.

\bigskip

\proclaim{\bf  5.1. The case $(E_6, \alpha_1)=(E_6, \alpha_6)$}
\endproclaim

The representation in question is $X = V(\omega_4 , D_5 )$, a half-spinor
representation for the group
${ G}=Spin (10)$. Here $\omega_4 = ({1\over 2}, {1\over 2}, {1\over 2}, {1\over 2}, {1\over 2})$.
$dim(X)=16$.
The weights of $X$ are vectors in 5 dimensional space, with coordinates equal to $\pm{1\over 2}$, with even number of negative coordinates.

We label the weight vectors in ${\goth g}_1$ by $[I]$ where $I$ is the subset of $\lbrace 1,2,3,4,5\rbrace$ of even cardinality where the sign of the component is negative.

The graded Lie algebra of type $E_6$ is
$${\goth g}(E_6)= {\goth g}_{-1}\oplus {\goth g}_0\oplus {\goth g}_1$$
with ${\goth g}_0= \CC\oplus {\goth so}(10)$.

The invariant scalar product $(,)$ on $\goth g$ restricted to ${\goth g}_1$ is given by the formula
$$([I], [J]) = 2-{1\over 2}\# (\lbrace I\setminus J\rbrace\cup\lbrace J\setminus I\rbrace) .$$

Thus possible scalar products are only $2,1,0$. So the possible root systems we can get are $A_1\times A_1$ and $A_1$.
Indeed, the triple product $A_1\times A_1\times A_1$ is not possible because there are no three subsets $I,J,K$ with cardinalities of three symmetric differences being 4.
This means we get

\proclaim {Proposition 5.1} The orbit structure is as follows.

$$\matrix number&diagram&dim&representative\\
0&0&0&0\\
1&A_1&11&[\emptyset]\\
2&2A_1&16&[\emptyset]+[\lbrace 1,2,3,4\rbrace]  \endmatrix$$

The containment diagram of corresponding orbit closures is
$$
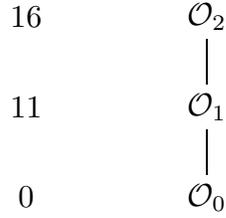
\xy
(15,0)*+{{\Cal O}_{0}}="o0";%
(15,10)*+{{\Cal O}_{1}}="o1";%
(15,20)*+{{\Cal O}_{2}}="o2";%
(-5,0)*{0};%
(-5,10)*{11};%
(-5,20)*{16};%
{\ar@{-} "o0"; "o1"};%
{\ar@{-} "o1"; "o2"};%
\endxy
$$
\endproclaim

The numerical data are as follows
$$\matrix number&degree&numerator\\
0&1&1\\
1&12&t^3+5t^2+5t+1\\
2&1&1
\endmatrix$$

The singularities data are as follows.

$$\matrix number&spherical&normal&C-M&R.S.&Gor\\
0&yes&yes&yes&yes&yes\\
1&yes&yes&yes&yes&yes\\
2&yes&yes&yes&yes&yes  \endmatrix$$

The coordinate ring of the closure $\overline{\OOO_1}$ of the orbit of the highest weight vector has the equivariant Betti table as follows

$$\matrix \CC&-&-&-&-&-\\
-&V_{\omega_1}&V_{\omega_5}&-&-&-\\
-&-&-&V_{\omega_4}&V_{\omega_1}&-\\
-&-&-&-&-&\CC
\endmatrix$$

Numerical Betti table is

$$\matrix 1&-&-&-&-&-\\
-&10&16&-&-&-\\
-&-&-&16&10&-\\
-&-&-&-&-&1
\endmatrix$$

The numerator of the Hilbert function is $t^3+5t^2+5t+1$.

The projectivization of the variety $\overline{\OOO_1}$ is self-dual. It is on the Knop-Mentzel list of highest weight vector orbits whose duals are degenerate.

\bigskip

\bigskip

\proclaim{\bf  5.2. The case $(E_6, \alpha_2)$}
\endproclaim

We have $X=\bigwedge^3 F$, $F=\CC^6$, ${ G}=GL(F)$.

The graded Lie algebra of type $E_6$ is
$${\goth g}(E_6)= {\goth g}_{-2}\oplus {\goth g}_{-1}\oplus {\goth g}_0\oplus {\goth g}_1\oplus {\goth g}_2$$
with ${\goth g}_0= \CC\oplus {\goth sl}(6)$,
${\goth g}_1=\bigwedge^3\CC^6$, ${\goth g}_2=\bigwedge^6\CC^6$.

\bigskip

The weights of ${\goth g}_1$ are $\epsilon_i +\epsilon_j+\epsilon_k$ for $1\le i<j<k\le 6$.
We label the corresponding weight vector by $[I]$ where $I$ is a cardinality 3 subset of $\lbrace 1,2,3,4,5,6\rbrace$.

The invariant scalar product on ${\goth g}$ restricted to ${\goth g}_1$ is
$$([I], [J])= \delta-1$$
where $\delta =\# (I\cap J)$.

We have four possible sets of weight vectors giving the roots of subalgebras $\goth s$.

\proclaim {Proposition 5.2} The representation ${\goth g}_1$ has five orbits.

$$\matrix number&diagram&dim&representative\\
0&0&0&0\\
1&A_1&10&[123]\\
2&2A_1&15&[123]+[145]\\
3&3A_1&19&[123]+[145]+[246]\\
4&A_2&20&[123]+[456]  \endmatrix$$

The containment diagram for the orbit closures is
$$
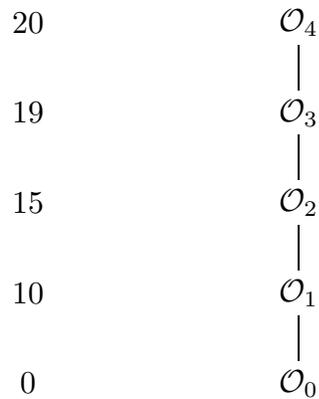
\xy
(15,0)*+{{\Cal O}_{0}}="o0";%
(15,10)*+{{\Cal O}_{1}}="o1";%
(15,20)*+{{\Cal O}_{2}}="o2";%
(15,30)*+{{\Cal O}_{3}}="o3";%
(15,40)*+{{\Cal O}_{4}}="o4";%
(-15,0)*{0};%
(-15,10)*{10};%
(-15,20)*{15};%
(-15,30)*{19};%
(-15,40)*{20};%
{\ar@{-} "o0"; "o1"};%
{\ar@{-} "o1"; "o2"};%
{\ar@{-} "o2"; "o3"};%
{\ar@{-} "o3"; "o4"};%
\endxy
$$

The geometric picture is as follows

$$\matrix number&projective\ picture&tensor\ picture\\
0&0&0\\
1&h.wt.vector&t=l_1\wedge l_2\wedge l_3\\
2&rank\ \le 5&t=l\wedge t'\\
3&\tau(\overline{\OOO_1})&&\\
4&\sigma(\overline{\OOO_1})&open  \endmatrix$$

The numerical data are as follows

$$\matrix number&degree&numerator\\
0&1&1\\
1&42&t^4+10t^3+20t^2+10t+1\\
2&42&t^5+5t^4+15t^3+15t^2+5t+1\\
3&4&t^3+t^2+t+1\\
4&1&1
\endmatrix$$

The algebraic pictures of singularities of orbit closures is as follows

$$\matrix number&spherical&normal&C-M&R.S.&Gor\\
0&yes&yes&yes&yes&yes\\
1&yes&yes&yes&yes&yes\\
2&yes&yes&yes&yes&yes\\
3&yes&yes&yes&yes&yes\\
4&no&yes&yes&yes&yes  \endmatrix$$
\endproclaim

We give the information on the minimal resolutions of the coordinate rings of the orbit closures.
Here $(a_1, a_2, a_3, a_4, a_5, a_6)$ denotes $S_{a_1, a_2, a_3, a_4, a_5, a_6} F^*$.

$\spadesuit$ The minimal resolution of the orbit closure $\overline{\OOO_1}$ (the cone over the Grassmannian $Grass(3,\CC^6)$) was first calculated in \cite{PW86}.
The transpose of its equivariant Betti table is
$$\matrix (0^6)&-&-&-&-\\
-&(2,1^4,0)&-&-&-\\
-&(3,2,1^4),(2^4,1,0)&-&-&-\\
-&(3^2,2^2,1^2)&(5,2^5),(3^5,0)&-&-\\
-&-&(5,3^3,2^2), (4^3,2^3),(4^2,3^3,1)&-&-\\
-&-&2*(5,4^2,3^2,2)&-&-\\
-&-&(5^2,4^3,2), (6,4^3,3^2), (5^3,3^3)&-&-\\
-&-&(5^5,2),(7,4^5)&(6^2,5^2,4^2)&-\\
-&-&-&(6^4,5,4), (7,6,5^4)&-\\
-&-&-&(7,6^4,5)&-\\
-&-&-&-&(7^6)
\endmatrix$$

The proof of the result is sketched in \cite{W03}, chapter 7, exercises 1-5. This coordinate ring is Gorenstein with rational singularities.

$\spadesuit$ The transpose of the equivariant Betti table of the coordinate ring of the orbit closure $\overline{\OOO_2}$ is
$$\matrix (0^6)&-&-&-&-&-\\
-&-&(2^3,1^3)&-&-&-\\
-&-&(3,2^4,1)&-&-&-\\
-&-&-&(4,3^4,2)&-&-\\
-&-&-&(4^3,3^3)&-&-\\
-&-&-&-&-&(5^6)\endmatrix$$

Numerical transpose of the Betti table is

$$\matrix 1&-&-&-&-&-\\
-&-&20&-&-&-\\
-&-&35&-&-&-\\
-&-&-&35&-&-\\
-&-&-&20&-&-\\
-&-&-&-&-&1\endmatrix$$

The coordinate ring $K[\overline{\OOO_2}]$ is Gorenstein with rational singularities.
The numerator of the Hilbert function is $t^5+5t^4+15t^3+15t^2+5t+1$.

$\spadesuit$ The coordinate ring of the orbit closure $\overline{\OOO_3}$ is a hypersurface of degree 4. It has rational singularities.

\proclaim {Remark 5.1}
This representation is one of the so-called subregular representations of Landsberg and Manivel \cite{LM01}.
The structure of the orbits and their defining ideals are uniform for these representations.
\endproclaim

\bigskip

\bigskip

\proclaim{\bf  5.3. The case $(E_6, \alpha_3)=(E_6,\alpha_5)$}
\endproclaim

We have $X=E\otimes\bigwedge^2 F$, $E=\CC^2$,
$F=\CC^5$,
${ G}= SL(E)\times SL(F)\times \CC^*$

The graded Lie algebra of type $E_6$ is
$${\goth g}(E_6)= {\goth g}_{-2}\oplus {\goth g}_{-1}\oplus {\goth g}_0\oplus {\goth g}_1\oplus {\goth g}_2$$
with ${\goth g}_0= \CC\oplus {\goth sl}(2)\oplus {\goth sl}(5)$,
${\goth g}_1=\CC^2\otimes\bigwedge^2\CC^5$, ${\goth g}_2=\bigwedge^2\CC^2\otimes\bigwedge^4\CC^5$.

Let $\lbrace e_1, e_2\rbrace$ be a basis of $E$, $\lbrace f_1 ,\ldots , f_5\rbrace$ be a basis of $F$.
We denote the tensor $e_a\otimes f_i\wedge f_j$ by $[a;ij]$.
The invariant scalar product on $\goth g$ restricted to ${\goth g}_1$ is
$$([a;ij], [b;kl])=\delta-1$$
where $\delta =\#(\lbrace a\rbrace\cap\lbrace b\rbrace )+\#(\lbrace i,j\rbrace\cap\lbrace k,l\rbrace ).$

The representation $X$ has eight  orbits.

$$\matrix number&{\goth s}&dim&representative\\
0&0&0&0\\
1&A_1&8&[1;12]\\
2&2A_1&11&[1;12]+[1;34] \\
3&2A_1&12&[1;12]+[2;13]\\
4&3A_1&15&[1;12]+[1;34]+[2;13]\\
5&A_2&16&[1;12]+[2;34]\\
6&A_2+A_1&18&[1;12]+[2;34]+[1;35].\\
7&A_2+2A_1&20&[1;12]+[2;34]+[1;35]+[2;15] \endmatrix$$

The containment diagram of the orbit closures is
$$
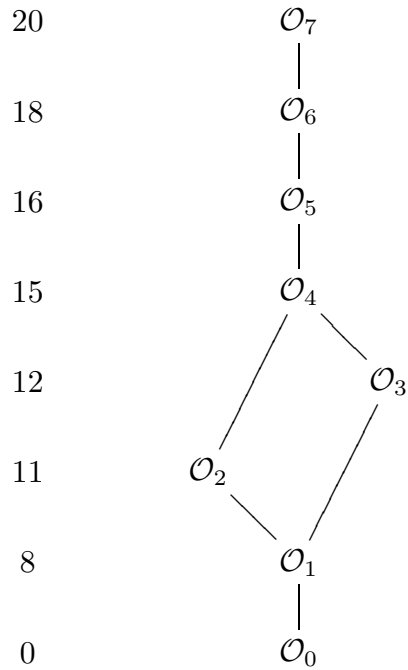
\xy
(15,0)*+{{\Cal O}_{0}}="o0";%
(15,10)*+{{\Cal O}_{1}}="o1";%
(5,20)*+{{\Cal O}_{2}}="o2";%
(25,30)*+{{\Cal O}_{3}}="o3";%
(15,40)*+{{\Cal O}_{4}}="o4";%
(15,50)*+{{\Cal O}_{5}}="o5";%
(15,60)*+{{\Cal O}_{6}}="o6";%
(15,70)*+{{\Cal O}_{7}}="o7";%
(-15,0)*{0};%
(-15,10)*{8};%
(-15,20)*{11};%
(-15,30)*{12};%
(-15,40)*{15};%
(-15,50)*{16};%
(-15,60)*{18};%
(-15,70)*{20};%
{\ar@{-} "o0"; "o1"};%
{\ar@{-} "o1"; "o2"};%
{\ar@{-} "o1"; "o3"};%
{\ar@{-} "o3"; "o4"};%
{\ar@{-} "o2"; "o4"};%
{\ar@{-} "o4"; "o5"};%
{\ar@{-} "o5"; "o6"};%
{\ar@{-} "o6"; "o7"};%
\endxy
$$

The geometry of the orbits is as follows.

$$\matrix number&proj.\ picture&tensor\  picture\\
0&0&0\\
1&C(Seg(\PP^1\times Grass(2,5)))&h.w.\ vector,&\\
&&t=l\otimes m_1\wedge m_2&\\
2&C(Seg(\PP^1\times P(\bigwedge^2 \CC^5)))&E-rank\le 1\\
3&general\ decomposable&t=l_1\otimes m\wedge m_1+l_2\otimes m\wedge m_2\\
4&\tau({\overline\OOO_1})&\\
5&\sigma_2 ({\overline\OOO}_1)&sum\ of\ 2\ decomposables,\\
&&F-rank\le 4&\\
6&&hyperdisc. \\
7&&generic \endmatrix$$

$$\matrix number&matrix\ picture&degeneration\\
1&rank\ 2, one\ var.&E,F-degenerate\\
2&one\ var.&E,F-degenerate\\
3&rank\ 2&F-degenerate\\
4&Pf_4 (X)= (x^2)&F-degenerate\\
5&Pf_4(X)=(xy)&F-degenerate\\
6&member\ of\ rank\ 2, \\
&Pf_4(X)=(x^2,xy)\\
7&Pf_4(X)=(x^2,xy,y^2)
\endmatrix$$

The numerical data are as follows

$$\matrix number&degree&numerator\\
0&1&1\\
1&35&1+12t+18t^2+4t^3 \\
2&10&1+9t\\
3&55&1+8t+  21t^2+  20t^3+  5t^4\\
4&64&1+5t+  15t^2+  25t^3+  15t^4+  3t^5 \\
5&30&5t^4+10t^3+10t^2+4t+1\\
6&15&4t^3+8t^2+2t+1\\
7&1&1
\endmatrix$$

The singularities data are as follows.

$$\matrix number&spherical&normal&C-M&R.S.&Gor\\
0&yes&yes&yes&yes&yes\\
1&yes&yes&yes&yes&no\\
2&yes&yes&yes&yes&no \\
3&yes&yes&yes&yes&no \\
4&yes&yes&yes&yes&no\\
5&no&yes&yes&yes&no\\
6&no&no&no&no&no\\
n(6)&no&yes&yes&yes&no\\
7&no&yes&yes&yes&yes \endmatrix$$

In the rest of this section we put $A=Sym(E^*\otimes\bigwedge^2 F^*)$ and we abbreviate $(a,b;c,d,e,f,g)$
for $S_{a,b}E^*\otimes S_{c,d,e,f,g}F^*$, and in resolutions this means the free $A$-module generated by this representation, i.e.
$S_{a,b}E^*\otimes S_{c,d,e,f,g}F^*\otimes A(-a-b)$.

Only the orbits $\OOO_7$ and $\OOO_6$ are non-degenerate.

\bigskip

$\spadesuit$ It is enough to deal with the orbit closure $\overline{\OOO_{6}}$ of codimension $2$.
This is the orbit closure given by pencils of antisymmetric
$5\times 5$ matrices with a member of rank $\le 2$.

As a desingularization  one can take the bundle  $\xi (6)$ over $\PP(E)\times Flag(1,3,4;C^5) $ spanned by weight vectors
$(1, 0;1, 1, 0, 0, 0)$, $(0, 1;1, 1, 0, 0, 0)$, $(1 ,0;1, 0 ,1, 0, 0)$,
\newline  $(0, 1;1, 0, 1, 0, 0)$, $(1, 0;0, 1, 1, 0, 0)$, $(1, 0;1, 0, 0, 1, 0)$,
$(0, 1;1, 0, 0, 1, 0)$ \newline and $ (1, 0;1, 0, 0, 0, 1)$.

To calculate the complex $\FF(6)_{\bullet}$, it is convinient to replace the bundle $\xi(6)$ with a bundle $\xi'(6)$ over
$\PP(E)\times Grass(2, F)$.
We treat each projective space as the Grassmannian with the tautological subbundles ${\Cal R}_E$, ${\Cal R}_F$
and tautological factor bundles ${\Cal Q}_E$ and ${\Cal Q}_F$ respectively.
The bundle $\xi' (6)$ is
$$\xi' (6)= {\Cal R}_E\otimes Ker(\bigwedge^2 F\rightarrow\bigwedge^2{\Cal Q}_F) + E\otimes\bigwedge^2{\Cal R}_F .$$

The complex $\FF(6)_{\bullet}$ is
$$0\rightarrow (4,1;2,2,2,2,2)\otimes A(-5)\rightarrow $$
$$\rightarrow (2,1;2,1,1,1,1)\otimes A(-3)\rightarrow (1,1;1,1,1,1,0)\otimes A(-2)\oplus  A.$$
The orbit closure is not normal.
The numerator of the Hilbert function is $4t^3+8t^2+2t+1$.
We have the exact sequence
$$0\rightarrow \CC[{\overline{\OOO_6}}]\rightarrow\CC[N({\overline{\OOO_6}})]\rightarrow C(6)\rightarrow 0.$$
In order to describe the defining ideal of $\overline{\OOO_6}$ we need to analyze the  $A$-module $C(6)$.
But considering the dual of the resolution of the coordinate ring of the orbit closure $\overline{\OOO_5}$ considered below, twisted by the appropriate power of the  determinant of $E$ (so this is in fact the resolution of the canonical module of the coordinate ring of $\overline{\OOO_5}$),
we get the acyclic complex $\FF (\bigwedge^2E\otimes\bigwedge^4{\Cal Q}_F)(5)_\bullet$ with the terms

$$0\rightarrow (5,5;4,4,4,4,4)\otimes A(-10)\rightarrow (4,3;3,3,3,3,2)\otimes A(-7)\rightarrow$$
$$\rightarrow (3,3;3,3,2,2,2)\otimes A(-6)\oplus (4,1;2,2,2,2,2)\otimes A(-5)\rightarrow$$
$$(2,1;2,1,1,1,1)\otimes A(-3)\rightarrow (1,1;1,1,1,1,0)\otimes A(-2).$$
It follows this complex is a free resolution of the $A$-module $C(6)$.
The resolution of the coordinate ring $\CC[{\overline \OOO_6}]$ is the mapping cone of a natural map
$$\FF(6)_\bullet\rightarrow \FF (\bigwedge^2E\otimes\bigwedge^4{\Cal Q}_F)(5)_\bullet$$
covering the natural surjection of $A$-modules.
There is a repeating representation $(4,1;2,2,2,2,2)$ occurring in both complexes. It has to cancel out. Indeed, if it does not, then it would
give part of minimal generators of the defining ideal. This is impossible, since $(4,1;2,2,2,2,2)$ does not occur in $S_5 (E\otimes\bigwedge^2 F)$.
We conclude that a minimal free resolution
of the $A$-module of $\CC[{\overline \OOO_6}]$ has the terms
$$0\rightarrow (5,5;4,4,4,4,4)\otimes A(-10)\rightarrow (4,3;3,3,3,3,2)\otimes A(-7)\rightarrow$$
$$\rightarrow  (3,3;3,3,2,2,2)\otimes A(-6)\rightarrow (0,0;0,0,0,0,0)\otimes A.$$
So generators come from (unique) irreducible representation $(3,3;3,3,2,2,2)$ occurring in
$S_6 (E\otimes\bigwedge^2 F)$.

\proclaim {Remark 5.2} The details of this calculation with the matrices in the complexes involved will appear in \cite{G11}.
\endproclaim

$\spadesuit$ Let us also comment on the generic $F$-degenerate  orbit $\OOO_{5}$ of codimension $4$.

The orbit closure is the set of the tensors of $F$-rank $\le 4$.
The bundle $\xi (5)$ is $E\otimes {\Cal Q}_F\otimes{\Cal R}_F$.
The complex $\FF(5)_\bullet$ has terms:
$$0\rightarrow (4,4;4,3,3,3,3)\otimes A(-8)\rightarrow (4,3;3,3,3,3,2)\otimes A(-7)\rightarrow$$
$$\rightarrow (4,1;2,2,2,2,2)\otimes A(-5)\oplus (2,2;2,2,2,1,1)\otimes A(-4)\rightarrow$$
$$\rightarrow (2,1;2,1,1,1,1\otimes A(-3))\rightarrow (0,0;0,0,0,0,0)\otimes A.$$

The orbit closure is normal with rational singularities.
The numerator of the Hilbert function is $5t^4+10t^3+10t^2+4t+1$.

\bigskip

\proclaim{\bf  5.4. The case $(E_6, \alpha_4)$}
\endproclaim

We have $X= E\otimes F\otimes H$, $E=\CC^2 ,F=H=\CC ^3$,
${G}=SL(E)\times SL(F)\times SL(H)\times \CC^*$.

The graded Lie algebra of type $E_6$ is
$${\goth g}(E_6)= {\goth g}_{-3}\oplus {\goth g}_{-2}\oplus {\goth g}_{-1}\oplus {\goth g}_0\oplus {\goth g}_1\oplus {\goth g}_2\oplus {\goth g}_3$$
with ${\goth g}_0= \CC\oplus {\goth sl}(2)\oplus{\goth sl}(3)\oplus{\goth sl}(3)$,
${\goth g}_1=\CC^2\otimes\CC^3\otimes\CC^3$, ${\goth g}_2=\bigwedge^2\CC^2\otimes\bigwedge^2\CC^3\otimes\bigwedge^2\CC^3$, ${\goth g}_3=S_{2,1}\CC^2\otimes\bigwedge^3\CC^3\otimes\bigwedge^3\CC^3$.

Let $\lbrace e_1, e_2\rbrace$ be a basis of $E$, and $\lbrace f_1, f_2, f_3\rbrace$,  $\lbrace h_1, h_2, h_3\rbrace$ bases of $F$, $H$ respectively.
We label $e_a\otimes f_i\otimes h_u$ by $[a;i;u]$.

The invariant scalar product on $\goth g$ restricted to ${\goth g}_1$ is
$$([a;i;u], [b;j;v])=\delta-1$$
where $\delta =\#(\lbrace a\rbrace\cap\lbrace b\rbrace )+\#(\lbrace i\rbrace\cap\lbrace j\rbrace )+\#(\lbrace u\rbrace\cap\lbrace v\rbrace ).$

The representation $X$ has 18 orbits.

$$\matrix number&\goth s &dim&representative\\
0&0&0&0\\
1&A_1&6&[1;1;1]\\
2&2A_1&8&[1;1;1]+[2;2;1]\\
3&2A_1&8&[1;1;1]+[2;1;2]\\
4&2A_1&9&[1;1;1]+[1;2;2]\\
5&3A_1&11&[1;1;1]+[1;2;2]+[2;1;2]\\
6&3A_1&10&[1;1;1]+[1;2;2]+[1;3;3]\\
7&A_2&12&[1;1;1]+[2;2;2]\\
8&A_2+A_1&13&[1;1;1]+[2;2;2]+[1;3;2]\\
9&A_2+A_1&13&[1;1;1]+[2;2;2]+[1;2;3] \\
10&A_2+2A_1&14&[1;1;1]+[2;2;2]+[1;3;2]+[2;3;1] \\
11&A_2+2A_1&14&[1;1;1]+[2;2;2]+[1;2;3]+[2;1;3] \\
12&A_2+2A_1&14&[1;1;1]+[2;2;2]+[1;2;3]+[1;3;2]\\
13&2A_2&14&[1;1;1]+[2;2;2]+[1;2;3]+[2;3;1]\\
14&A_3&15&[1;1;1]+[2;2;2]+[1;3;3]\\
15&2A_2+A_1&16&[1;1;1]+[2;2;2]+[1;2;3]+[2;3;1]+[1;3;2]\\
16&A_3+A_1&17&[1;1;1]+[2;2;2]+[1;3;3]+[2;1;3]\\
17&D_4(a_1)&18&[1;1;1]+[2;1;1]+[1;2;2]-[2;2;2]+[1;3;3]  \endmatrix$$

Let us comment on the Lie algebra $\goth s$ of type $D_4(a_1)$. This is a Lie algebra of type $D_4$ with the grading taking values
on simple roots as follows: $deg(\alpha_1 )=deg(\alpha_3 )=deg(\alpha_4 )=1, deg(\alpha_2 )=0$.
This algebra has a decomposition
$${\goth s}={\goth s}_{-3}\oplus{\goth s}_{-2}\oplus {\goth s}_{-1}\oplus {\goth s}_0\oplus {\goth s}_1\oplus{\goth s}_2\oplus{\goth s}_3$$
with $dim\ {\goth g}_0=dim\ {\goth g}_1 =6$, $dim\ {\goth g}_2=3$, $dim\ {\goth g}_3 =3$. Thus the six roots in degree $1$ have
certain pattern of scalar products. There is one choice (up to permutations by Weyl group of $G_0$) of such six roots in ${\goth g}_1$.
They are $[1;1;1]$, $[2;1;1]$, $[1;2;2]$, $[2;2;2]$, $[1;3;3]$, $[2;3;3]$. Taking a generic element in the subspace spanned by these roots we get
a point in the generic orbit.

Our representation has two interpretations: It is the space of $3\times 3$ matrices of linear forms in two variables (matrix picture),
or it is the space of representations of dimension vector $(3,3)$ of the Kronecker quiver (quiver picture). In particular, in the quiver picture the
$SL(F)\times SL(H)$-nullcone of our space is a complete intersection given by the condition $det=0$. It has three irreducible components which are orbit closures
$\overline{\OOO_9}, \overline{\OOO_{10}}, \overline{\OOO_{12}}$.

$$\matrix number&proj.\ picture&matrix\ pic.&quiver\ pic.\\
0&0&0&0\\
1&C(Seg (\PP^1\times \PP^2\times \PP^2 ))&rank\ 1,one\ var.&\\
2&&$H$-rank\le 1&\\
3&&$F$-rank\le 1&\\
4&&$E$-rank\le 1,$F$-$H$-rank\le 2&\\
5&\tau(\overline\OOO_1 )&\\
6&&$E$-rank\le 1&\\
7&\sigma_2 (\overline\OOO_1 )&\\
8&hyperdet\ in\ \overline{\OOO_{10}}&\\
9&hyperdet\ in\  \overline{ \OOO_{11}}& \\
10&&$H$-rank\le 2&(3,2) subrep. \\
11&&$F$-rank\le 2&(1,0) subrep. \\
12&&rank-1-member,det=x^3&\\
13&&&(2,1) subrep.\\
14&&rank-1-member,det=x^2y&\\
15&&det=x^3&\\
16&&det=x^2y&\\
17&\sigma_3 (\overline\OOO_1 )&generic&  \endmatrix$$

The containment diagram for orbit closures is

$$
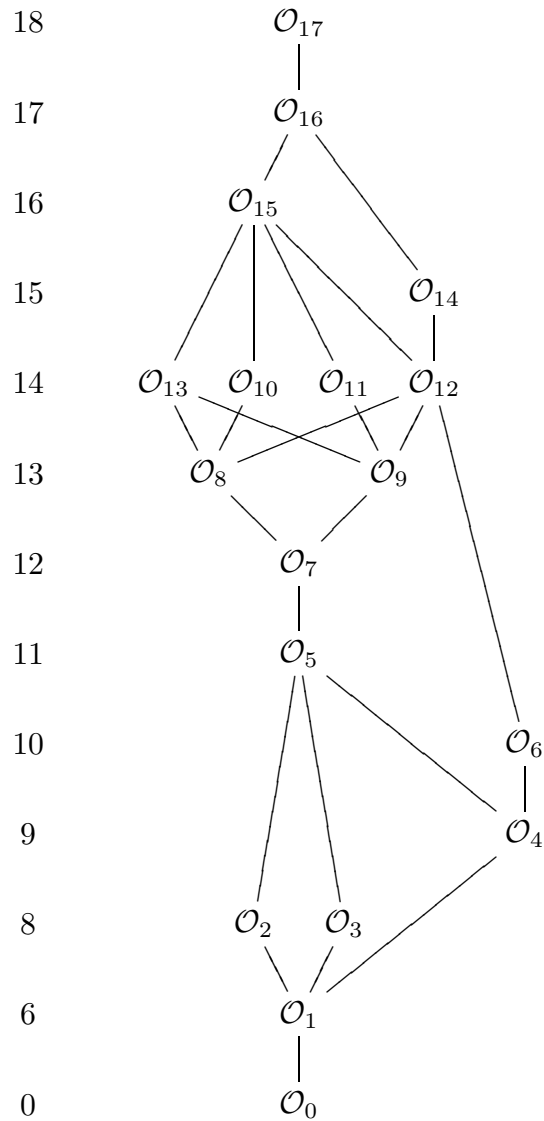
\xy
(15,0)*+{{\Cal O}_{0}}="o0";%
(15,10)*+{{\Cal O}_{1}}="o1";%
(10,20)*+{{\Cal O}_2}="o2";%
(20,20)*+{{\Cal O}_3}="o3";%
(40,30)*+{{\Cal O}_4}="o4";%
(15,50)*+{{\Cal O}_5}="o5";%
(40,40)*+{{\Cal O}_6}="o6";%
(15,60)*+{{\Cal O}_7}="o7";%
(5,70)*+{{\Cal O}_8}="o8";%
(25,70)*+{{\Cal O}_9}="o9";%
(10,80)*+{{\Cal O}_{10}}="o10";%
(20,80)*+{{\Cal O}_{11}}="o11";%
(30,80)*+{{\Cal O}_{12}}="o12";%
(0,80)*+{{\Cal O}_{13}}="o13";%
(30,90)*+{{\Cal O}_{14}}="o14";%
(10,100)*+{{\Cal O}_{15}}="o15";%
(15,110)*+{{\Cal O}_{16}}="o16";%
(15,120)*+{{\Cal O}_{17}}="o17";%
(-15,0)*{0};%
(-15,10)*{6};%
(-15,20)*{8};%
(-15,30)*{9};%
(-15,40)*{10};%
(-15,50)*{11};%
(-15,60)*{12};%
(-15,70)*{13};%
(-15,80)*{14};%
(-15,90)*{15};%
(-15,100)*{16};%
(-15,110)*{17};%
(-15,120)*{18};%
{\ar@{-} "o0"; "o1"};%
{\ar@{-} "o1"; "o2"};%
{\ar@{-} "o3"; "o1"};%
{\ar@{-} "o4"; "o1"};%
{\ar@{-} "o2"; "o5"};%
{\ar@{-} "o3"; "o5"};%
{\ar@{-} "o4"; "o5"};%
{\ar@{-} "o4"; "o6"};%
{\ar@{-} "o5"; "o7"};%
{\ar@{-} "o6"; "o12"};%
{\ar@{-} "o7"; "o8"};%
{\ar@{-} "o7"; "o9"};%
{\ar@{-} "o8"; "o12"};%
{\ar@{-} "o8"; "o13"};%
{\ar@{-} "o8"; "o10"};%
{\ar@{-} "o9"; "o12"};%
{\ar@{-} "o9"; "o13"};%
{\ar@{-} "o9"; "o11"};%
{\ar@{-} "o12"; "o14"};%
{\ar@{-} "o12"; "o15"};%
{\ar@{-} "o13"; "o15"};%
{\ar@{-} "o10"; "o15"};%
{\ar@{-} "o11"; "o15"};%
{\ar@{-} "o14"; "o16"};%
{\ar@{-} "o15"; "o16"};%
{\ar@{-} "o16"; "o17"};%
\endxy
$$

This is seen directly from the geometric descriptions and from knowing the representatives.

\proclaim{Remark 5.3} The partial order above is displayed in the paper by Parfenov \cite{P01}.
However in his table one of the adjacencies is not marked (in his diagram he needs the
edge connecting vertices $9$ and $11$). This can be seen by noticing that his adjacency diagram
is not invariant under the involution exchanging two factors $\CC^3$.
\endproclaim

The numerical data is as follows

$$\matrix number&degree&numerator\\
0&1&1\\
1&30&1+12t+15t^2+2t^3\\
2&21&1+10t+10t^2\\
3&21&1+10t+10t^2\\
4&24&1+9t+9t^2+5t^3\\
5&108&1+7t+28t^2+48t^3+21t^4+3t^5\\
6&9&1+8t\\
7&57&1+6t+21t^2+20t^3+9t^4\\
8&72&1+5t+24t^2+26t^3+16t^4\\
9&72&1+5t+24t^2+26t^3+16t^4\\
10&15&1+4t+10t^2\\
11&15&1+4t+10t^2\\
12&66&1+4t+19t^2+22t^3+16t^4+4t^5\\
13&51&1+4t+19t^2+18t^3+9t^4\\
14&24&1+3t+15t^2+5t^3\\
15&27&1+2t+3t^2+4t^3+5t^4+6t^5+4t^6+2t^7\\
16&12&1+t+t^2+3t^3+3t^4+3t^5\\
17&1&1
\endmatrix$$

The data regarding singularities is

$$\matrix number&spherical&normal&C-M&R.S.&Gor\\
0&yes&yes&yes&yes&yes\\
1&yes&yes&yes&yes&no\\
2&yes&yes&yes&yes&no\\
3&yes&yes&yes&yes&no\\
4&yes&yes&yes&yes&no\\
5&yes&yes&yes&yes&no\\
6&yes&yes&yes&yes&no\\
7&no&yes&yes&yes&no\\
8&no&yes&yes&yes&yes\\
9&no&yes&yes&yes&yes \\
10&no&yes&yes&yes&no \\
11&no&yes&yes&yes&no \\
12&no&no&no&no&no\\
n(12)&no&yes&yes&yes&no\\
13&no&no&no&no&no\\
n(13)&no&yes&yes&yes&no\\
14&no&no&no&no&no\\
n(14)&no&yes&yes&yes&no\\
15&no&yes&yes&yes&no\\
16&no&no&yes&no&yes\\
n(16)&no&yes&yes&yes&yes\\
17&no&yes&yes&yes&yes  \endmatrix$$

Only the last six orbits are $F$-non-degenerate and $H$-non-degenerate.
Pairs of orbits of the same codimension are symmetric with respect to the involution switching $F$ and $H$.
The determinantal orbits, the highest weight vector orbit, the  zero orbit and the generic orbit are easy to handle.

Below we describe the desingularizations of the non-degenerate orbit closures and related complexes $\FF_\bullet$
one can construct from them using geometric technique from \cite{W03}.

We will abbreviate $(a,b;c,d,e;f,g;h)$ for $S_{a,b}E^*\otimes S_{c,d,e} F^*\otimes S_{f,g,h}H^*$.
All our complexes $\FF_\bullet$ are written over the ring $A=Sym(E^*\otimes F^*\otimes H^*)$.
In order to relate these complexes to each other we will use the notation $\FF(\Cal V)(s)_\bullet$ to denote the twisted complex
$\FF(\Cal V )_\bullet$ related to our chosen desingularization of the orbit closure $\overline{\OOO(s)}$ ($\Cal V$ is a twisting vector bundle) .

\bigskip

$\spadesuit$ The hyperdeterminant orbit $\OOO_{16}$.

This is the hypersurface given by the $2\times 3\times 3$ matrices with vanishing hyperdeterminant. Its desingularization lives on $\PP(E)\times \PP(F)\times \PP(G)$.
We treat each projective space as the set of lines with the tautological subbundles ${\Cal R}_E$, ${\Cal R}_F$ and ${\Cal R}_H$ respectively.
The bundle $\xi$ is
$$\xi = {\Cal R}_E\otimes  {\Cal R}_F\otimes H+  {\Cal R}_E\otimes F\otimes {\Cal R}_H + E\otimes  {\Cal R}_F\otimes {\Cal R}_H .$$

The complex $\FF(16)_{\bullet}$ is
$$0\rightarrow (4,2;2,2,2;2,2,2)\otimes A(-6)\rightarrow (2,1;1,1,1;1,1,1)\otimes A(-3)\oplus A.$$
The determinant of this matrix is the hyperdeterminant $\Delta$ of our matrix which has degree $12$.

\bigskip

$\spadesuit$  The codimension $2$ orbit $\OOO_{15}$.

The orbit closure has a geometric description as the set of pencils of $3\times 3$ matrices whose determinant is a cube of a linear form.

The bundle $\xi$ is the bundle over $\PP^1\times Flag(1,2;\CC^3)\times Flag(1,2,\CC^3) $ corresponding to  submodule with the weights
$$\matrix (0,1;0,0,1;0,0,1), & (0,1;0,1,0;0,0,1), & (0,1;0,0,1;0,1,0), \\
(0,1;0,1,0;0,1,0), & (0,1;1,0,0;0,0,1), & (0,1;0,0,1;1,0,0),\\
(1,0;0,0,1;0,0,1), & (1,0;0,1,0;0,0,1), & (1,0;0,0,1;0,1,0).
\endmatrix $$
This bundle can be thought of as the set of weights in the graphic form as follows.

$$ \matrix X&X&X\\X&X&O\\X&O&O\endmatrix\ \ \ \ \  \matrix X&X&O\\X&O&O\\O&O&O\endmatrix$$
Here the first matrix represents the weights with  $(1,0)$  on the first coordinate, and the second matrix represents the weghts with$(0,1)$
on the first coordinate. The symbol $X$ denotes the weight in $\eta$, the symbol $O$ -- the weight in $\xi$.
Our incidence space $Z(15)$ has dimension $9+3+3+1=16$, so it projects on the orbit of codimension $2$.

The calculation of cohomology reveals that we get a complex $\FF(15)_\bullet$

$$0\rightarrow (5,4;3,3,3;3,3,3)\otimes A(-9)\rightarrow (4,2;2,2,2;2,2,2)\otimes A(-6)\rightarrow A.$$
This complex can be obtained also from the hyperdeterminant complex by looking at the kernel of the transpose of the cubic part
of the complex $\FF(16)_\bullet$. We see the complex is determinantal so it must give a resolution of the reduced ideal. This is in fact a proof that the incidence space $Z(15)$ is a desingularization.
Dividing by the regular sequence of $16$ generic linear forms gives a resolution of the ring with Hilbert function
$${1-3t^6+2t^9\over (1-t)^2}=1+2t+3t^2+4t^3+5t^4+6t^5+4t^6+2t^7 .$$
This means the degree of the orbit closure is $21$.

\bigskip

$\heartsuit$ The codimension $3$ orbit $\OOO_{14}$.

This orbit closure has a nice geometric description. It consists of pencils of $3\times 3$ matrices of linear forms containing a matrix of rank $\le 1$.

The bundle $\xi$ is the submodule with the weights
$$\matrix
(0,1;0,0,1;0,0,1), &  (0,1;0,1,0;0,0,1), & (0,1;1,0,0;0,0,1), \\
(0,1;0,0,1;0,1,0), & (0,1;0,1,0;0,1,0), & (0,1;1,0,0;0,1,0), \\
(0,1;0,0,1;1,0,0), &  (0,1;0,1,0;1,0,0), &  (0,1;1,0,0;1,0,0), \\
(1,0;0,0,1;0,0,1).&& \endmatrix $$
This bundle can be thought of as the set of weights in the graphic form as follows.

$$ \matrix X&X&X\\X&X&X\\X&X&O\endmatrix \ \ \ \ \ \matrix O&O&O\\O&O&O\\O&O&O\endmatrix$$
Here the first matrix represents the weights with $(1,0)$  on the first coordinate, and the second matrix represents the weghts with $(0,1)$
on the first coordinate. The symbol $X$ denotes the weight in $\eta$, the symbol $O$ -- the weight in  $\xi$.
Our incidence space has dimension $10+2+2+1=15$, so it projects on the orbit of codimension $3$. It is clearly
$\overline{\OOO_{14}}$.

The calculation of cohomology reveals that we get a complex $\FF(14)_\bullet$

$$0\rightarrow (5,1;2,2,2;2,2,2)\otimes A(-6)\rightarrow (3,1;2,1,1;2,1,1)\otimes A(-4)\rightarrow $$
$$\rightarrow(2,1;1,1,1;2,1,0)\otimes A(-3)\oplus (2,1;2,1,0;1,1,1)\otimes A(-3)\rightarrow $$
$$\rightarrow (1,1;1,1,0;1,1,0)\otimes A(-2)\oplus A.$$

The terms of the complex $F(14)_\bullet$ can be understood by looking at the partitions on two last coordinates. They give terms of the Gulliksen-Negard complex resolving $2\times 2$ minors of the generic $3\times 3$ matrix. The corresponding term on the first coordinate comes from the symmetric power on the bundle ${\Cal R}_E$. This means the term comes from $H^1$ on the $\PP(E)$ coordinate (except the trivial term), hence the shift in homological degree.
The orbit is not normal.

The cokernel was resolved using Macaulay 2 (see \cite{G11}). Its Betti table is
$$\matrix 9&32&27&.&.&.&.&.&.\\
.&.&.&5&.&.&.&.&.\\
.&.&104&342&477&372&180&54&7\\
.&.&.&.&.&.&.&.&.\\
.&.&.&.&.&.&1&.&.&
\endmatrix
   $$
   From this it is not difficult to see that the resolution of the coordinate ring of $\CC[{\overline\OOO_{14}}]$ has a Betti table
   $$\matrix 1&.&.&.&.&.&.&.\\
   .&.&.&.&.&.&.&.\\
   .&.&.&.&.&.&.&.\\
   .&.&.&.&.&.&.&.\\
.&.&.&.&.&.&.&.\\
.&104&342&477&372&180&54&7\\
.&.&.&.&.&.&.&.\\
.&.&.&.&.&1&.&.&
\endmatrix
   $$

   This implies that $\CC[\overline{\OOO_{14}}]$ is not Cohen-Macaulay.

   Let us describe the cokernel $C(14)$ more precisely. This is a very interesting case, as this is the first case of a non-normal orbit closure whose non-normality locus $\overline{\OOO_{12}}$ is also non-normal.

   Consider the epimorphism
   $$N(\CC[\overline{\OOO_{14}}])\rightarrow N(\CC [\overline{\OOO_{14}}]).$$
   It extends to the exact sequence
   $$0\rightarrow K(14,12)\rightarrow N(\CC[\overline{\OOO_{14}}])\rightarrow N(\CC [\overline{\OOO_{14}}])\rightarrow 0.\eqno{(*)}$$

   It induces the epimorphism of the cokernels $C(14)\rightarrow C(12)$.
   The cokernel $C(12)$ and its resolution is described below in the analysis of the orbit closure $\overline{\OOO_{12}}$.
   We have an exact sequence
   $$0\rightarrow L(14,12)\rightarrow C(14)\rightarrow C(12)\rightarrow 0.$$
   The kernel $L(14,12)$ turns out to be a twisted module supported on the orbit closure $\overline{\OOO_6}$ with the twist $(2,1;1,1,1;1,1,1)$.
   Both modules $C(12)$ and $L(14,12)$ are maximal Cohen-Macaulay modules.

   In fact there are two desingularizations $Z'_{14}$ and $Z'_{12}$
   of $\overline{\OOO_{14}}$ and of $\overline{\OOO_{12}}$ which are compatible, i.e. the bundle $\eta'(12)$ is a factor of $\eta '(14)$, with the kernel of
   rank one. There are in fact two such compatible pairs. Thus the exact sequence $(*)$ is just the sequence of sections of the relative Koszul complex in one variable.

   The module $K(14,12)$ is a twisted module on $Z'(14)$ and it is a maximal Cohen-Macaulay module supported on this  orbit closure.

   All of this leads to the exact sequence
   $$0\rightarrow J(14,12)\rightarrow K(14,12)\rightarrow L(14,12)\rightarrow 0$$
   where $J(14,12)$ is the defining ideal of $\overline{\OOO_{12}}$ in $\CC[{\overline{\OOO_{14}}}]$.

   Let us describe the generators of the defining ideal of $\CC[\overline{\OOO_{14}}]$. Looking at the covariant $S_2 E^*\otimes\bigwedge^2 F^*\otimes\bigwedge^2 H^*$ we see (after substituting the representative of
   our orbit $\OOO_{14}$ that we get $9$ binary quadratic forms with a common root.

   Looking at the coordinate ring $A'=Sym (S_2E^*\otimes W^*)$ where $W=\bigwedge^3 F\otimes\bigwedge^3 H$ and applying geometric technique for the subvariety of several quadrics with a common root
   ($\xi =S_2\RRR_E\otimes W$, $\eta = \QQQ_E\otimes \RRR_E\otimes W$),
   we get another nonnormal stratum. We can resolve the cokernel there and we find out that the stratum has the defining ideal generated by cubics, namely by representation
   $S_{3,3}E^*\otimes\bigwedge^3 W^*$.  This gives $84$ equations.
   The other equations come from representations $(3,3;2,2,2;4,1,1)$ and $(3,3;4,1,1;2,2,2)$ which occur with multiplicity 1 in $A_6$.

\bigskip

$\spadesuit$ The codimension $4$ orbit $\OOO_{13}$.

The geometric description of this orbit can be understood from another point of view. Our representation can be thought of
as the set of representations of Kronecker quiver of dimension vector $(3,3)$.

$$\CC^3\ \matrix \rightarrow\\ \rightarrow \endmatrix \ \CC^3 .$$

Our orbit closure then is the set of representations having a subrepresentation of dimension vector $(2,1)$.
The desingularization $Z(13)$ lives on $Grass(2,F)\times Grass(2,H)$ and the bundle $\xi$ is

$$\xi = E\otimes{\Cal R}_F\otimes{\Cal R}_H .$$
Graphically we have

$$ \matrix X&X&X\\X&O&O\\X&O&O\endmatrix\ \ \ \ \  \matrix X&X&X\\X&O&O\\X&O&O\endmatrix .$$

Our incidence space has dimension $10+2+2=14$, so it projects on the orbit of codimension $4$.
In this case one can see directly that $Z(13)$ is a desingularization of $\overline{\OOO_{13}}$.

The calculation of Euler characteristics suggests  that we get a complex $\FF(13)_\bullet$

$$0\rightarrow (4,4;3,3,2;3,3,2)\otimes A(-8)\rightarrow (4,3;3,2,2;3,2,2)\otimes A(-7)\rightarrow $$
$$\rightarrow(2,2;2,1,1;2,1,1)\otimes A(-4)\oplus (3,1;2,1,1;2,1,1)\otimes A(-4)\oplus $$
$$\oplus (3,2;2,2,1;2,2,1)\otimes A(-5)\rightarrow $$
$$\rightarrow (2,1;1,1,1;1,1,1)\otimes A(-3)\oplus (2,1;1,1,1;2,1,0)\otimes A(-3)\oplus $$
$$\oplus (2,1;2,1,0;1,1,1)\otimes A(-3)\oplus (3,0;1,1,1;1,1,1)\otimes A(-3) \rightarrow$$
$$\rightarrow (1,1;1,1,0;1,1,0))\otimes A(-2)\oplus (0,0;0,0,0;0,0,0)\otimes A.$$

The orbit closure is obviously not normal.

The complex $\FF(13)_\bullet$ gives a minimal resolution of the normalization
$N(\overline{\OOO_{13}})$. We have an exact sequence
$$0\rightarrow \CC[{\overline{\OOO_{13}}}]\rightarrow \CC[N({\overline{\OOO_{13}}})]\rightarrow C(13)\rightarrow 0.$$
The complex $\FF(13)_\bullet$ reveals that the $A$-module $C(13)$ has the presentation
$$(2,1;1,1,1;1,1,1)\otimes A(-3)\oplus(2,1;1,1,1;2,1,0)\otimes A(-3)\oplus $$
$$\oplus (2,1;2,1,0;1,1,1)\otimes A(-3)\rightarrow (1,1;1,1,0;1,1,0)\otimes A(-2).$$
Indeed, the representation $(3,0;1,1,1;1,1,1)$ can map only into the trivial term in the complex $\FF(13)_\bullet$.

We can look for the module with the above presentation among twisted modules supported in smaller orbits.
It turns out the right choice is the orbit closure $\overline{\OOO_7}$. Its bundle $\xi$ has a diagram
$$ \matrix X&X&O\\X&X&O\\O&O&O\endmatrix\ \ \ \ \  \matrix X&X&O\\X&X&O\\O&O&O\endmatrix .$$
This bundle lives on the space $Grass(1, F)\times Grass(1, H)$. Consider the twisted complex
$\FF (\bigwedge^2 E\otimes \bigwedge^2{\Cal Q}_F\otimes\bigwedge^2{\Cal Q}_H )(7)_\bullet$.
Its terms are:

$$0$$
$$\downarrow$$
$$ (6,6;4,4,4;4,4,4)$$
$$\downarrow$$
$$(5,4;3,3,3;4,3,2)\oplus(5,4;4,3,2;3,3,3)\oplus (6,3;3,3,3;3,3,3)$$
$$\downarrow$$
$$(4,4;3,3,2;3,3,2)\oplus (4,4;3,3,2;4,2,2)\oplus (4,4;4,2,2;3,3,2)\oplus 2*(5,3;3,3,2;3,3,2)$$
$$\downarrow$$
$$(4,3;3,2,2;3,2,2)\oplus (4,3;3,2,2;3,3,1)\oplus (4,3;3,3,1;3,2,2)\oplus (5,1;2,2,2;2,2,2)$$
$$\downarrow$$
$$(2,2;2,1,1;2,1,1)\oplus (3,1;2,1,1;2,1,1)\oplus (3,3;2,2,2;3,3,0)\oplus (3,3;3,3,0;2,2,2)$$
$$\downarrow$$
$$(2,1;1,1,1;1,1,1)\oplus (2,1;1,1,1;2,1,0)\oplus(2,1;2,1,0;1,1,1)$$
$$\downarrow$$
$$(1,1;1,1,0;1,1,0).$$

The resolution of he $A$-module $\CC[{\overline{\OOO_{13}}}]$ can be constructed as a mapping cone of the map
$$\FF(13)_\bullet\rightarrow \FF (\bigwedge^2 E\otimes \bigwedge^2{\Cal Q}_F\otimes\bigwedge^2{\Cal Q}_H )(7)_\bullet$$
covering the natural epimorphism of $A$-modules.
The mapping cone is not a minimal resolution but the repeating representations might  cancel out.
Let us see that the pairs of representations $(2,2;2,1,1;2,1,1)$ and $(3,1;2,1,1;2,1,1)$ cancel out.
From this we can deduce that the defining ideal of $\overline{\OOO_{13}}$ is generated by the representations
$(3,0;1,1,1;1,1,1)$ in degree $3$ and the representations $(3,3;2,2,2;3,3,0)$ and $(3,3;3,3,0;2,2,2)$ in degree $6$.

Indeed, if $(3,1;2,1,1;2,1,1)$ would not cancel out, it would contribute to the minimal generators of the defining ideal. However
it occurs once in $S_4 (E\otimes F\otimes H)$ so that representation is already in the ideal generated by $(3,0;1,1,1;1,1,1)$.
Regarding representation $(2,2;2,1,1;2,1,1)$, if it would occur in the defining ideal of $\overline{\OOO_{13}}$ the analysis of the next section will show that then the
defining ideal of $\overline{\OOO_{13}}$ would contain the defining ideal of $\overline{\OOO_{12}}$. But this is impossible since both orbit closures have the same dimension.
The homological dimension of $\CC[{\overline{\OOO_{13}}}]$ as an $A$-module equals $5$ because the top of the resolution of $C(13)$
does not cancel out. So this coordinate ring is not Cohen-Macaulay. Also the non-normality locus of $\overline{\OOO_{13}}$ equals to $\overline{\OOO_7}$.

\bigskip

$\spadesuit$ The codimension $4$ orbit $\OOO_{12}$.

The geometric description of this orbit closure is that these are pencils of matrices with a matrix of rank $\le 1$,
and with determinant which is a cube of a linear form.

The desingularization $Z(12)$ is given by the bundle $\xi$ which in graphical form is

$$ \matrix X&X&X\\X&X&O\\X&O&O\endmatrix\ \ \ \ \  \matrix X&O&O\\O&O&O\\O&O&O\endmatrix .$$

Our incidence space $Z(12)$ has dimension $7+3+3+1=14$, so it projects on an orbit closure of codimension $4$.
It is clearly a different orbit than the other orbits of codimension 4 as the determinant of a general matrix in the pencil  is non-zero.
In this case one can see directly that $Z(12)$ is a desingularization of $\overline{\OOO_{12}}$.

The calculation of cohomology reveals that we get a complex $\FF(12)_\bullet$

$$0\rightarrow (6,3;3,3,3;3,3,3)\rightarrow (4,3;3,2,2;3,2,2)\oplus (5,1;2,2,2;2,2,2)\rightarrow $$
$$\rightarrow(3,1;2,1,1;2,1,1)\oplus (3,3;2,2,2;3,2,1)\oplus (3,3;3,2,1;2,2,2)\rightarrow $$
$$\rightarrow (2,1;1,1,1;1,1,1)\oplus (2,1;1,1,1;2,1,0)\oplus (2,1;2,1,0;1,1,1) $$
$$\rightarrow (1,1;1,1,0;1,1,0))\oplus (0,0;0,0,0;0,0,0).$$

The orbit closure is obviously not normal. The complex $\FF(12)_\bullet$ gives a minimal resolution of the normalization
$N(\overline{\OOO_{12}})$. We have an exact sequence
$$0\rightarrow \CC[{\overline{\OOO_{12}}}]\rightarrow \CC[N({\overline{\OOO_{12}}})]\rightarrow C(12)\rightarrow 0.$$
It follows that the cokernel $C(12)$ is the same as the cokernel $C(13)$ considered in the previous section.
We deduce that the resolution of he $A$-module $\CC[{\overline{\OOO_{12}}}]$ can be constructed as a mapping cone of the map
$$\FF(12)_\bullet\rightarrow \FF (\bigwedge^2 E\otimes \bigwedge^2{\Cal Q}_F\otimes\bigwedge^2{\Cal Q}_H )(7)_\bullet$$
covering the natural epimorphism of $A$-modules.
The mapping cone is not a minimal resolution but the repeating representations might  cancel out.
The only doubtful pair is $(3,1;2,1,1;2,1,1)$. But one can see directly that this representation does not vanish on $\overline{\OOO_{12}}$.
So this pair cancels out and the generators of the defining ideal of the orbit closure $\overline{\OOO_{12}}$ are representations
$(2,2;2,1,1;2,1,1)$ in degree $4$ and the representations $(3,3;2,2,2;3,3,0)$ and $(3,3;3,3,0;2,2,2)$ in degree $6$.
The homological dimension of $\CC[{\overline{\OOO_{12}}}]$ as an $A$-module equals $5$ because the top of the resolution of $C(12)$
does not cancel out. So this coordinate ring is not Cohen-Macaulay. Also the non-normality locus of $\overline{\OOO_{12}}$ equals to $\overline{\OOO_7}$.

\bigskip\bigskip

\head \S 6. The type I $F_4$ representations. \endhead

\bigskip\bigskip

\proclaim{\bf  6.1. The case $(F_4, \alpha_1)$}
\endproclaim

We have  $X= V(\omega_3 ,C_3 )$, ${ G}= SP(F)$ with $F$ a
symplectic space of dimension 6.

The graded Lie algebra of type $F_4$ is
$${\goth g}(F_4)= {\goth g}_{-2}\oplus {\goth g}_{-1}\oplus {\goth g}_0\oplus {\goth g}_1\oplus {\goth g}_2$$
with ${\goth g}_0= \CC\oplus {\goth sp}(6)$,
${\goth g}_1=V(\omega_3 ,C_3)$, ${\goth g}_2=\CC$.

The roots  in ${\goth g}_1$  correspond to weights of $V(\omega_3 ,C_3)$ as follows:
$$(1,0,0,0)\leftrightarrow\epsilon_1+\epsilon_2+\epsilon_3 ,$$
$$(1,1,0,0)\leftrightarrow \epsilon_1+\epsilon_2-\epsilon_3 ,$$
$$(1,1,1,0)\leftrightarrow \epsilon_1-\epsilon_2+\epsilon_3 ,$$
$$(1,1,1,1)\leftrightarrow \epsilon_1 ,\ (1,1,2,0)\leftrightarrow -\epsilon_1+\epsilon_2+\epsilon_3 ,$$
$$(1,1,2,1)\leftrightarrow \epsilon_2 ,\ (1,2,2,0)\leftrightarrow -\epsilon_3 ,$$
$$(1,1,2,2)\leftrightarrow \epsilon_3 ,\ (1,2,2,1)\leftrightarrow -\epsilon_2 ,$$
$$(1,2,2,2)\leftrightarrow -\epsilon_1 ,\ (1,2,3,1)\leftrightarrow \epsilon_1-\epsilon_2-\epsilon_3 ,$$
$$(1,2,3,2)\leftrightarrow -\epsilon_1+\epsilon_2-\epsilon_3 ,$$
$$(1,2,4,2)\leftrightarrow -\epsilon_1-\epsilon_2+\epsilon_3 ,$$
$$(1,3,4,2)\leftrightarrow -\epsilon_1-\epsilon_2-\epsilon_3 .$$

There are 5 orbits.

$$\matrix number&\goth s&representative\\
0&0&0\\
1&A_1&(1,0,0,0)\\
2&{\tilde A}_1&(1,1,1,1)\\
3&A_1+{\tilde A}_1&(1,0,0,0)+(1,2,3,1)\\
4&A_2&(1,0,0,0)+(1,3,4,2)
\endmatrix$$

The containment diagram of  orbit closures is linear
$$
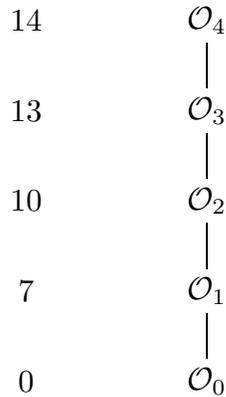
\xy
(15,0)*+{{\Cal O}_{0}}="o0";%
(15,10)*+{{\Cal O}_{1}}="o1";%
(15,20)*+{{\Cal O}_{2}}="o2";%
(15,30)*+{{\Cal O}_{3}}="o3";%
(15,40)*+{{\Cal O}_{4}}="o4";%
(-5,0)*{0};%
(-5,10)*{7};%
(-5,20)*{10};%
(-5,30)*{13};%
(-5,40)*{14};%
{\ar@{-} "o0"; "o1"};%
{\ar@{-} "o1"; "o2"};%
{\ar@{-} "o2"; "o3"};%
{\ar@{-} "o3"; "o4"};%
\endxy
$$

In terms of the symplectic group the $\goth g$ invariant scalar product is (???).

This is another term in subexceptional series of Landsberg and Manivel.
There is one invariant $\Delta$ of degree 4 and the orbits are as follows

$$\matrix number&dim&description&geometry\\
0&0&0&\\
1&7&h.weight\ vector&\\
2&10&Sing(\nabla )&\\
3&13&\nabla&\\
4&14&generic&
  \endmatrix$$

  The numerical data are as follows

$$\matrix number&degree&numerator\\
0&1&1\\
1&16&t^3+7t^2+7t+1\\
2&21&6t^3+10t^2+4t+1\\
3&4&t^3+t^2+t+1\\
4&1&1
\endmatrix$$

The algebraic pictures of singularities of orbit closures is as follows

$$\matrix number&spherical&normal&C-M&R.S.&Gor\\
0&yes&yes&yes&yes&yes\\
1&yes&yes&yes&yes&yes\\
2&yes&yes&yes&yes&no\\
3&yes&yes&yes&yes&yes\\
4&no&yes&yes&yes&yes  \endmatrix$$

\bigskip

$\spadesuit$ The codimension $1$ orbit $\OOO_{3}$.
 The orbit closure $\overline{\OOO_3}=Y^\vee_{hw}$ is normal, the defining ideal is generated by  $\Delta $.

 \bigskip

$\spadesuit$ The codimension $4$ orbit $\OOO_{2}$.
The orbit closure $\overline{\OOO_2}=Sing (Y^\vee_{hw})$  whose defining ideal is generated by a covariant
of degree 3 consisting of partial derivatives of $\Delta$. The minimal free resolution of the defining ideal is
$$0\rightarrow V_{\omega_1}\otimes A(-7)\rightarrow V_{\omega_2}\otimes A(-6)\rightarrow $$
$$\rightarrow V_{2\omega_1}\otimes A(-4)\rightarrow V_{\omega_3}\otimes A(-3)\rightarrow A.$$

\bigskip

$\spadesuit$ The codimension $7$ orbit $\OOO_{1}$.
 The orbit $\overline{\OOO_1}=Y_{hw}$ of the highest weight vector. It is
an orbit of codimension  7.
The minimal resolution of the defining ideal is

$$0\rightarrow A(-10)\rightarrow V_{2\omega_1}\otimes A(-8)\rightarrow V_{\omega_1+\omega_2}\otimes A(-7)\rightarrow V_{\omega_1+\omega_3}\otimes A(-6)\rightarrow$$
$$\rightarrow V_{\omega_1+\omega_3}\otimes A(-4)\rightarrow V_{\omega_1+\omega_2}\otimes A(-3)\rightarrow V_{2\omega_1}\otimes A(-2)\rightarrow A.$$

\bigskip

\proclaim{\bf  6.2. The case $(F_4, \alpha_2)$}
\endproclaim

We have    $X=E\otimes S_2 F$ with $E=\CC^2
,F=\CC^3$, ${G}=SL(E)\times SL(F)\times \CC^*$.

The graded Lie algebra of type $F_4$ is
$${\goth g}(F_4)= {\goth g}_{-3}\oplus {\goth g}_{-2}\oplus {\goth g}_{-1}\oplus {\goth g}_0\oplus {\goth g}_1\oplus {\goth g}_2\oplus {\goth g}_3$$
with ${\goth g}_0= \CC\oplus {\goth sl}(2)\oplus{\goth sl}(3)$,
${\goth g}_1=\CC^2\otimes S_2 \CC^3$, ${\goth g}_2=\bigwedge^2\CC^2\otimes S_{2,2} \CC^3$, ${\goth g}_3=S_{2,1}\CC^2\otimes S_{2,2,2}\CC^3$.

$$\matrix number&\goth s&representative\\
0&0&0\\
1&A_1&(0,1,0,0)\\
2&{\tilde A}_1&(0,1,1,0)\\
3&A_1+{\tilde A}_1&(0,1,0,0)+(0,1,2,1)\\
4&A_1+{\tilde A}_1& (0,1,0,0)+(1,1,1,0)\\
5&A_2&(0,1,0,0)+(1,1,2,0)\\
6&{\tilde A}_2&(0,1,2,1)+(1,1,1,1)\\
7&A_2+{\tilde A}_1&(0,1,0,0)+(1,1,2,0)+(0,1,2,1)\\
8&B_2&(0,1,0,0)+(1,1,2,1)\\
9&{\tilde A}_2+A_1&(1,1,1,1,)+(0,1,2,1)+(0,1,0,0)\\
10&C_3(a_1)&<(0,1,0,0), (0,1,1,0), (1,1,1,1), (0,1,2,0), (1,1,2,1)>\\
11&F_4(a_3)&(0,1,0,0)+(1,1,0,0)+(0,1,2,0)-(1,1,2,0)+(0,1,2,2)
\endmatrix$$

\proclaim{Remark} The support algebra of the nilpotent in the orbit $O_9$
has simple roots $(1,1,2,1), -(0,0,1,0)$ and $(0,1,2,0)$.
\endproclaim

The containment diagram of orbit closures is
$$
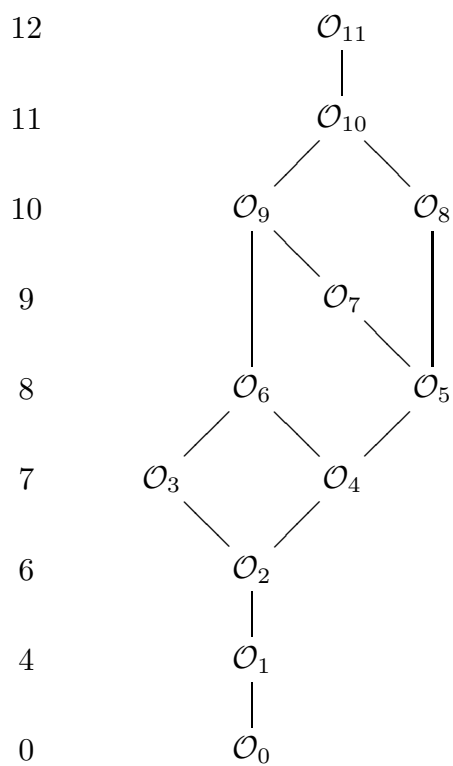
\xy
(15,0)*+{{\Cal O}_{0}}="o0";%
(15,10)*+{{\Cal O}_{1}}="o1";%
(15,20)*+{{\Cal O}_{2}}="o2";%
(5,30)*+{{\Cal O}_{3}}="o3";%
(25,30)*+{{\Cal O}_{4}}="o4";%
(35,40)*+{{\Cal O}_{5}}="o5";%
(15,40)*+{{\Cal O}_{6}}="o6";%
(25,50)*+{{\Cal O}_{7}}="o7";%
(35,60)*+{{\Cal O}_{8}}="o8";%
(15,60)*+{{\Cal O}_{9}}="o9";%
(25,70)*+{{\Cal O}_{10}}="o10";%
(25,80)*+{{\Cal O}_{11}}="o11";%
(-10,0)*{0};%
(-10,10)*{4};%
(-10,20)*{6};%
(-10,30)*{7};%
(-10,40)*{8};%
(-10,50)*{9};%
(-10,60)*{10};%
(-10,70)*{11};%
(-10,80)*{12};%
{\ar@{-} "o0"; "o1"};%
{\ar@{-} "o1"; "o2"};%
{\ar@{-} "o2"; "o3"};%
{\ar@{-} "o2"; "o4"};%
{\ar@{-} "o3"; "o6"};%
{\ar@{-} "o4"; "o6"};%
{\ar@{-} "o4"; "o5"};%
{\ar@{-} "o5"; "o7"};%
{\ar@{-} "o5"; "o8"};%
{\ar@{-} "o6"; "o9"};%
{\ar@{-} "o7"; "o9"};%
{\ar@{-} "o8"; "o10"};
{\ar@{-} "o9"; "o10"};%
{\ar@{-} "o10"; "o11"};%
\endxy
$$

$$\matrix number&dim&description&geometry\\
0&0&0&\\
1&4&h.w.\ vector&1\ variable,\ rank\ 1\\
2&6&&1\ variable,\ rank\ 2\\
3&7&&1\ variable,\ rank\ 3\\
4&7&&pencil\ of\ 2\times 2\ matrices\ with\ 0\ discriminant\\
5&8&&pencil\ of\ 2\times 2\ matrices\\
6&8&&pencil\ of\ matrices\ with\ 2\times 2\ block\ of\ zeros\\
7&9&&rank\ 1\ member, det=x^3\\
8&10&&rank\ 1\ member, det=x^2y\\
9&10&&det=x^3,\\
10&11&\nabla&det=x^2y\\
11&12&generic&
  \endmatrix$$

  The numerical data are:

  $$\matrix number&degree&numerator\\
  0&1&1\\
  1&12&1+8t+ 3t^2\\
  2&15&1+6t+6t^2+2t^3\\
  3&6&1+5t\\
  4&36&1+11t+18t^2+6t^3\\
  5&15&1+4t+10t^2\\
  6&21&1+10t+10t^2\\
  7&32&1+3t+ 12t^2+  10t^3+  6t^4\\
  8&12&1+2t+9t^2\\
  9&27&1+2t+6t^2+10t^3+8t^4\\
  10&12&1+t+4t^2+6t^3\\
  11&1&1
  \endmatrix$$

  The singularities data are

  $$\matrix number&spherical&normal&C-M&Rat.\ sing.&Gorenstein\\
  0&yes&yes&yes&yes&yes\\
  1&yes&yes&yes&yes&no\\
  2&yes&yes&yes&yes&no\\
  3&yes&yes&yes&yes&no\\
  4&yes&no&yes&no&no\\
  5&no&yes&yes&yes&no\\
  6&no&no&no&no&no\\
  n(6)&no&yes&yes&yes&no\\
  7&no&no&no&no&no\\
  n(7)&no&yes&yes&yes&no\\
  8&no&no&no&no&no\\
  n(8)&no&yes&yes&yes&no\\
  9&no&no&yes&no&no\\
  n(9)&no&yes&yes&yes&no\\
  10&no&no&yes&no&yes\\
  n(10)&no&yes&yes&yes&no\\
  11&yes&yes&yes&yes&yes
  \endmatrix$$

Notice that the Lie algebra $\goth s$ of type $C_3(a_1)$ has the grading 1 on roots $\alpha_1 ,\alpha_3$ of $C_3$ and grading zero on $\alpha_2$.
The grading is
$${\goth s}= {\goth s}_{-3}\oplus {\goth s}_{-2}\oplus {\goth s}_{-1}\oplus {\goth s}_0\oplus {\goth s}_1\oplus {\goth s}_2\oplus {\goth s}_3$$
with the dimensions $dim\ {\goth s}_0 =dim\ {\goth g}_1=5$, $dim\ {\goth s}_2 =2$, $dim\ {\goth s}_3=1$.
The algebra $F_4 (a_3)$ is just our graded algebra $\goth g$ so by definition the representative is a generic element of ${\goth g}_1$.
It can be for example $(0,1,0,0)+(1,1,0,0)+(0,1,2,0)-(1,1,2,0)+(0,1,2,2)$.

The roots correspond to the weight vectors in ${\goth g}_1$ as follows
$$(0,1,0,0)\leftrightarrow e_1\otimes f_1^2 ,$$
$$(1,1,0,0)\leftrightarrow e_2\otimes f_1^2 , (0,1,1,0)\leftrightarrow e_1\otimes f_1f_2 ,$$
$$(1,1,1,0)\leftrightarrow e_2\otimes f_1f_2 , (0,1,1,1)\leftrightarrow e_1\otimes f_1f_3 , (0,1,2,0)\leftrightarrow e_1\otimes f_2^2 ,$$
$$(1,1,1,1)\leftrightarrow e_2\otimes f_1f_3 ,(1,1,2,0)\leftrightarrow e_2\otimes f_2^2 , (0,1,2,1)\leftrightarrow e_1\otimes f_2f_3 ,$$
$$(1,1,2,1)\leftrightarrow e_2\otimes f_2f_3 ,(0,1,2,2)\leftrightarrow e_1\otimes f_3^2 ,$$
$$(1,1,2,2)\leftrightarrow e_2\otimes f_3^2 .$$

Here is more concrete description of the orbits. We treat the elements of our representation
as $3\times 3$ symmetric matrices whose entries are the homogeneous  linear forms in the variables $x, y$.
Alternatively we think of them as pencils of $3\times 3$ symmetric matrices.

We give some information on the minimal free resolutions of the coordinate rings of some orbit closures.
We set $A= Sym (E^*\otimes S_2 F^*)$ and we abbreviate $(a_1, a_2;b_1, b_2, b_3)$ for
$S_{a_1, a_2}E^*\otimes S_{b_1, b_2, b_3}F^*$.

 \bigskip

$\spadesuit$ The codimension $1$ orbit $\OOO_{10}$.

 The orbit $\overline{\OOO_9}=Y^\vee_{hw}$ closure. It has codimension 1, but it's not normal.
The invariant is the discriminant of the determinant of our matrix, it has degree 12.
The closure of the orbit are the pencils whose determinant has a double root.
The orbit closure is not normal, the minimal resolution of the coordinate ring of the normalization is
$$0\rightarrow (2,2;4,2,2)\otimes A(-4)\rightarrow (2,1;2,2,2)\otimes A(-3)\oplus (1,1;2,1,1)\otimes A(-2)\oplus A.$$

 \bigskip

$\spadesuit$ The codimension $2$ orbit $\OOO_{9}$. The orbit of codimension 2 of the pencils whose determinant is a cube of linear form,
but they have no member of rank 1. The resolution of the normalization of the coordinate ring is
$$0\rightarrow (3,3;5,4,3)\otimes A(-6)\rightarrow (3,2;4,3,3)\otimes A(-5)\oplus (2,2;4,2,2)\otimes A(-4)\rightarrow$$
$$\rightarrow (1,1;2,1,1)\otimes A(-2)\oplus A.$$

The resolution of the coordinate ring of the orbit closure is
$$0\rightarrow (5,4;6^3)\otimes A(-6)\rightarrow (4,2;4^3)\otimes A(-4)\rightarrow (0^2;0^3)\otimes A.$$
It is the resolution of the rational normal curve for the binary cubic when we substitute for the binary cubic the covarant $(3,0;2^3)$.

 \bigskip

$\spadesuit$ The codimension $2$ orbit $\OOO_{8}$.
 The orbit closure $\overline{\OOO_8}$ of codimension 2 consisting of pencils with
one member of rank $\le 1$.
Here the determinant of the pencil has a double root.
The resolution of the normalization of the coordinate ring is
$$0\rightarrow (3,1;3,3,2)\otimes A(-4)\rightarrow (2,1;3,2,1)\otimes A(-3)\rightarrow (1,1;2,2,0)\otimes A(-2)\oplus A.$$

 \bigskip

$\spadesuit$ The codimension $3$ orbit $\OOO_{7}$.
Here determinant of the pencil is the power of linear form and there is a member of rank 1.
The resolution of the normalization of the coordinate ring is
$$0\rightarrow (4,3;5,5,4)\otimes A(-7)\rightarrow (3,3;5,4,3)\otimes A(-6)\oplus (3,1;3,3,2)\otimes A(-4)\rightarrow$$
$$\rightarrow (2,1;3,2,1)\otimes A(-3)\oplus (2,1;2,2,2)\otimes A(-3)\rightarrow (1,1;2,2,0)\otimes A(-2)\oplus A.$$

 \bigskip

$\spadesuit$ The codimension $4$ orbit $\OOO_{6}$.
 The orbit closure $\overline{\OOO_6}$ consists  of pencils which (after change of variables) have a $2\times 2$ block of zeros in the upper left corner,
 The resolution of the normalization of the coordinate ring is
 $$0\rightarrow (3,3;5,5,2)\otimes A(-6)\rightarrow (3,2;5,3,2)\otimes A(-5)\rightarrow$$
 $$\rightarrow (3,1;4,2,2)\otimes A(-4)\oplus (2,2;3,3,2)\otimes A(-4)\oplus (2,1;4,1,1)\otimes A(-3)\rightarrow$$
 $$\rightarrow (3,0;2,2,2)\otimes A(-3)\oplus (1,1;3,1,0)\otimes A(-2)\oplus (2,0;2,1,1)\otimes A(-2)\rightarrow $$
 $$\rightarrow (1,0;1,1,0)\otimes A(-1)\oplus A.$$

\bigskip

$\spadesuit$ The codimension $4$ orbit $\OOO_{5}$.
The orbit closure $\overline{\OOO_5}$ consists of pencils which (after change of variables) are the pencils of $2\times 2$ symmetric matrices,
The orbit closure is normal, the resolution of the coordinate ring is
$$0\rightarrow (3,3;6,3,3)\otimes A(-6)\rightarrow (3,2;4,3,3)\otimes A(-5)\oplus (3,2;5,3,2)\otimes A(-5)\rightarrow$$
$$\rightarrow (3,1;4,2,2)\otimes A(-4)\oplus (2,2;4,3,1)\otimes A(-4)\oplus $$
$$\oplus  (2,2;3,3,2)\otimes A(-4)\oplus (3,1;3,3,2)\otimes A(-4)\rightarrow $$
$$\rightarrow (3,0;2,2,2)\otimes A(-3)\oplus (2,1;3,2,1)\otimes A(-3)\oplus A.$$

\bigskip

$\spadesuit$ The codimension $5$ orbit $\OOO_{4}$.

\bigskip

$\spadesuit$ The codimension $5$ orbit $\OOO_{3}$.
  This is the  generic $E$-degenerate orbit  of codimension 5.

\bigskip

$\spadesuit$ The codimension $6$ orbit $\OOO_{2}$.
This is the  codimension 6 $E$-degenerate orbit of matrices of rank 2.

\bigskip

$\spadesuit$ The codimension $8$ orbit $\OOO_{1}$.
 This is the  highest weight orbit $Y_{hw}$, it has codimension 8.

\bigskip

\proclaim{\bf  6.3. The case $(F_4, \alpha_3)$}
\endproclaim

We have   $X=E\otimes  F$ with $E=\CC^2
,F=\CC^3$, ${G}=SL(E)\times SL(F)\times \CC^*$.

The graded Lie algebra of type $F_4$ is
$${\goth g}(F_4)= {\goth g}_{-4}\oplus {\goth g}_{-3}\oplus {\goth g}_{-2}\oplus {\goth g}_{-1}\oplus {\goth g}_0\oplus {\goth g}_1\oplus {\goth g}_2\oplus {\goth g}_3\oplus {\goth g}_4$$
with ${\goth g}_0= \CC\oplus {\goth sl}(2)\oplus{\goth sl}(3)$,
${\goth g}_1=\CC^2\otimes  \CC^3$, ${\goth g}_2= S_2 \CC^2\otimes \bigwedge^2 \CC^3$, ${\goth g}_3=S_{2,1}\CC^2\otimes \bigwedge^3\CC^3$, ${\goth g}_4= S_{2,2}\CC^2\otimes S_{2,1,1}\CC^3$.
The orbit closures in this representations are just the determinantal varieties, treated for example in \cite{W03}, chapter 6.

\bigskip

\proclaim{\bf  6.4. The case $(F_4, \alpha_4)$}
\endproclaim

We have   $X= V(\omega_3 ,B_3 )$, ${ G}= Spin (7)$.

The graded Lie algebra of type $F_4$ is
$${\goth g}(F_4)= {\goth g}_{-2}\oplus {\goth g}_{-1}\oplus {\goth g}_0\oplus {\goth g}_1\oplus {\goth g}_2$$
with ${\goth g}_0= \CC\oplus {\goth so}(7)$,
${\goth g}_1=V(\omega_3 , B_3)$, ${\goth g}_2=\CC$.

The correspondance of roots in ${\goth g}_1$ with the weight vectors in $V(\omega_3 , B_3 )$ is as follows.

$$(0,0,0,1)\leftrightarrow {1\over 2}(\epsilon_1+\epsilon_2+\epsilon_3 ),$$
$$(0,0,1,1)\leftrightarrow {1\over 2}(\epsilon_1+\epsilon_2-\epsilon_3 ),$$
$$(0,1,1,1)\leftrightarrow {1\over 2}(\epsilon_1-\epsilon_2+\epsilon_3 ),$$
$$(1,1,1,1)\leftrightarrow {1\over 2}-(\epsilon_1+\epsilon_2+\epsilon_3 ), (0,1,2,1)\leftrightarrow {1\over 2}-(\epsilon_1-\epsilon_2-\epsilon_3 )$$
$$(1,1,2,1)\leftrightarrow {1\over 2}(-\epsilon_1+\epsilon_2-\epsilon_3 ),$$
$$(1,2,2,1)\leftrightarrow {1\over 2}(-\epsilon_1-\epsilon_2+\epsilon_3 ),$$
$$(1,2,3,1)\leftrightarrow {1\over 2}(-\epsilon_1-\epsilon_2-\epsilon_3 ).$$

All roots in ${\goth g}_1$ are short, so there are only three orbits

$$\matrix number&\goth s&representative\\
0&0&0\\
1&{\tilde A}_1&(0,0,0,1)\\
2&{\tilde A}_2&(0,0,0,1)+(1,2,3,1)
\endmatrix$$

$$\matrix number&dim&description&geometry\\
0&0&0&\\
1&7&h.weight\ vector=\nabla&\\
2&8&generic&
  \endmatrix$$

\bigskip

$\spadesuit$ The codimension $1$ orbit $\OOO_{1}$.
The highest weight orbit closure $\overline{\OOO_1}$ has codimension 1, its defining ideal is generated by an invariant of degree 2,

\bigskip\bigskip

\head \S 7. The type I  $G_2$ representations. \endhead

\bigskip\bigskip

\proclaim{\bf  7.1. The case $(G_2, \alpha_1)$}
\endproclaim

We have $X=  E$, $E=\CC^2$, ${ G}=GL(E)$.

The graded Lie algebra of type $G_2$ is
$${\goth g}(G_2)= {\goth g}_{-3}\oplus {\goth g}_{-2}\oplus {\goth g}_{-1}\oplus {\goth g}_0\oplus {\goth g}_1\oplus {\goth g}_2\oplus{\goth g}_3$$
with ${\goth g}_0= \CC\oplus {\goth sl}(2)$,
${\goth g}_1=  \CC^2$, ${\goth g}_2=\bigwedge^2\CC^2, {\goth g}_3=S_{2,1}\CC^2$.

The orbits in this representation are the generic one and the origin.

\bigskip

\proclaim{\bf  7.2. The case $(G_2, \alpha_2)$}
\endproclaim

We have $X= S_3 E$, $E=\CC^2$, ${ G}=GL(E)$.

The graded Lie algebra of type $G_2$ is
$${\goth g}(G_2)= {\goth g}_{-2}\oplus {\goth g}_{-1}\oplus {\goth g}_0\oplus {\goth g}_1\oplus {\goth g}_2$$
with ${\goth g}_0= \CC\oplus {\goth sl}(2)$,
${\goth g}_1= S_3 \CC^2$, ${\goth g}_2=\CC$.

The correspondance of roots in ${\goth g}_1$ with the weight vectors in $S_3 E$ is as follows.
$$(0,1)\leftrightarrow e_1^3 ,$$
$$(1,1)\leftrightarrow e_1^2e_2 ,$$
$$(2,1)\leftrightarrow e_1e_2^2 ,$$
$$(3,1)\leftrightarrow e_2^3 .$$

There are four orbits

$$\matrix number&\goth s&representative\\
0&0&0\\
1&A_1&(0,1)\\
2&{\tilde A}_1&(1,1)\\
3&G_2(a_1 )&(0,1)+(3,1)
\endmatrix$$

$$\matrix number&dim&description&geometry\\
0&0&0&\\
1&2&h.weight\ vector&\\
2&3&\nabla&\\
3&4&generic&
  \endmatrix$$

The ring of $SL(E)$-invariants is generated by the discriminant $\Delta$ of degree 4.

We give some information on the minimal free resolutions of the coordinate rings of some orbit closures.
We set $A= Sym (S_3 E^*)$ and we abbreviate $(a_1, a_2)$ for
$S_{a_1, a_2}E^*$.

\bigskip

$\spadesuit$ The codimension $1$ orbit $\OOO_{2}$.

 The orbit closure $\overline{\OOO_2}=Y^\vee_{hw}$  of cubics with vanishing discriminant. It is not normal. The resolution of its normalization is
 $$0\rightarrow (4,2)\otimes A[-2] \rightarrow (2,1)\otimes A[-1] \oplus (0,0)\otimes A[0].$$

\bigskip

$\spadesuit$ The codimension $2$ orbit $\OOO_{1}$.The orbit closure $\overline{\OOO_1}=Y_{hw}$ of the orbit of the highest weight vector, it has codimension 2.
It is a rational normal curve, it has rational singularities and the resolution of the coordinate ring is the Eagon-Northcott complex
$$0\rightarrow (5,4)\otimes A(-3)\rightarrow (4,2)\otimes A(-2)\rightarrow A.$$

\bigskip\bigskip

\head \S 7. Conclusions. \endhead

Here are some general conclusions we checked type by type.

Let $X_n$ be the Dynkin diagram of type $E_6, F_4, G_2$. The orbit  $X_{disc.}$ has the following properties.

\proclaim {Proposition 7.1}
\item {a)} The variety $X_{disc.}$ is closed $G_0$-equivariant and irreducible and thus it is an orbit closure,
\item{b)} The variety $X_{disc.}$ is a hypersurface if and only if the ring of invariants \break $Sym({\goth g}_1^*)^{(G,G)}$ contains a non-constant invariant. In these cases we have
$$Sym({\goth g}_1^*)^{(G,G)}=K[\Delta ]$$
 and $X_{disc.}$ is a hypersurface given by vanishing of the invariant $\Delta$,
\item{c)} In the cases when $Sym({\goth g}_1^*)^{(G,G)}=K$ the variety $X_{disc.}$ has codimension bigger than one.
\endproclaim

\proclaim {Proposition 7.2}
Let $X_n$ be of type $E_6, F_4, G_2$.
\item{a)} The orbit closure is spherical if and only if the support algebra $\goth s$ has all simple components of type $A_1$,
\item{b)} If the orbit closure $\overline{\OOO_v}$ is not normal then its normalization is contained in the representation ${\goth g}_1\oplus {\goth g}_i$ for some $i>1$,
where ${\goth g}_i$ is the $i$-th graded component in the grading associate to the simple root $\alpha_k$.
\endproclaim


\newpage

\enddocument